\documentclass[lettersize,journal]{IEEEtran}
\usepackage{amsmath,amsfonts}
\usepackage{array}
\usepackage[caption=false,font=normalsize,labelfont=sf,textfont=sf]{subfig}
\usepackage{bbm}
\usepackage{textcomp}
\usepackage{stfloats}
\usepackage{url}
\usepackage{verbatim}
\usepackage{graphicx}
\usepackage{cite}
\usepackage{booktabs}
\usepackage{array}
\usepackage{pifont}
\usepackage{tabularx}
\usepackage{multirow}
\usepackage{algorithm}
\usepackage{algorithmic}
\usepackage{color}
 
\makeatletter
\newcommand{\removelatexerror}{\let\@latex@error\@gobble}
\makeatother

\newtheorem{theorem}{Theorem}
\newtheorem{lemma}{Lemma}
\newtheorem{assumption}{Assumption}
\newtheorem{corollary}{Corollary}
\newtheorem{example}{Example}[section]

\begin{document}

\title{Efficient Algorithm Design of Dynamic Spectrum Access by Whittle Index}

\author{Keqin Liu$^*$, Qizhen Jia, Yiying Zhang, and Zhi Ding 
\thanks{The first and second authors are with School of Mathematics and Physics, Xi'an Jiaotong-Liverpool University, Suzhou, China. The third author is with ByteDance, Shanghai, China. The fourth author is with Department of Electrical and Computer Engineering, University of California at Davis, USA.}
\thanks{Part of this work was presented at 2023 International Conference on Statistics, Applied Mathematics and Computing Science. This manuscript is also available at arXiv:2501.00236.}
\thanks{$^*$Corresponding author, keqin.liu@xjtlu.edu.cn}
}



\maketitle

\begin{abstract}
Dynamic spectrum access problem is an important problem that allows a wireless sub-network to use channels temporarily unoccupied by the parent network for minimizing the spectrum waste. Previous work has shown that the sequential channel allocation problem for the sub-network can be formulated within the restless multi-armed bandits (RMAB) framework. The objective is to maximize the expected long-term return over an infinite horizon while minimizing interference to the parent network. Different from the previous work that exploits a binary feedback (e.g., ACK/NAK) to compensate for sensing errors, we leverage the finer and more robust channel quality indicator (CQI) feedback to update the information state (belief vector) of the sub-network. However, the implementation of CQI-based observation model yields significantly more complex belief transition behaviors in an infinite state space and worsens the curse of dimensionality of dynamic programming. To overcome this challenge, we dive into the rich structures of the value functions and obtain tight bounds on their derivatives. These results lead to the proof of optimality of threshold policies on a single-channel problem with subsidy and subsequently a closed-form channel index function using an iterative method to approximate the well-known Whittle index policy, which offers a low-complexity solution for ranking the currently available channels whose states are never directly observable. Through extensive numerical studies, we demonstrate the superior performance and robustness of our proposed algorithm.

\end{abstract}


\section{Introduction}
\subsection{Dynamic Spectrum Access}
\IEEEPARstart{D}{ynamic} spectrum access (DSA) refers to the process of efficiently allocating radio resources within a sub-network to ensure optimal performance and quality of service (QoS) for its users. For example, in 4G-HeNB (4G-home eNodeB) network, a small, low-power cellular base station (like a femtocell) provides localized coverage within a limited area, typically in urban high-density areas, by sharing a set of wireless channels with the main/parent network. Effective sequential resource scheduling involves the coordination and allocation of available channels (temporarily not used by the main network) to multiple user devices connected to the small base station~\cite{3GPP_HeNB}. The main challenge here is that the sub-network cannot perfectly observe the availabilities of chosen channels and the channel availabilities are themselves time-varying. In order to achieve high efficiency in radio resource sharing for future 5G deployment of DSA, the primary objective is to maximize sub-network throughput, minimize interference with the parent network, and enhance overall network performance over long-run~\cite{GTI_FemtoWP}.

To achieve efficient channel allocation for DSA, various algorithms and techniques are employed in different network structures~\cite{anand_machine_2023, lin_interference_2020, scopelliti_energy-quality_2020, pu_novel_2020, eslami_joint_2023, famili_ofdra_2023, sharma_secrecy_2023, zhang_joint_2022, xu_robust_2022, li_distributed_2022}. These include channel allocation and interference management. Channel allocation requires that the wireless scheduler assigns appropriate channels to user devices based on factors such as channel quality, interference levels, and user priority. The allocation is dynamic, depending on the network requirements and traffic conditions. Meanwhile, the scheduler also needs to mitigate interference to other networks sharing the same set of channels. Interference coordination techniques are employed to optimize spectral efficiency and minimize cross-network interference.

Overall, DSA plays a crucial role in optimizing resource utilization, managing interference, and providing a satisfactory user experience within the localized coverage area of the sub-network. It enables efficient and reliable communication for a diverse range of user devices while maintaining network performance objectives. In this paper, we mainly focus on the optimization of channel allocation for maximizing the expected data throughput of the sub-network over long-run.

\subsection{Restless Multi-armed Bandit Problem}
In this paper, we will formulate our channel allocation problem as a restless multi-armed bandit process (RMAB). RMAB is a generalization of the classical multi-armed bandit problem (MAB). MAB is a well-known mathematical model that serves as a foundational framework for dynamic resource allocation problems and has been widely used in multiple fields \cite{odeyomi_online_2023, sneh_radar_2023, shisher_learning_2023}. In the classical formulation, a player is challenged by the task of selecting a single arm out of $N$ options, subsequently receiving a random reward dependent on the state of the arm. The chosen arm undergoes a state transition according to a Markovian rule while other arms do not change their states. At all times, the states of all arms are perfectly observable. The player's goal is to maximize their cumulative discounted reward over an infinite time interval according to a specific arm selection policy. The inception of the general MAB concept can be traced back to its original exposition in 1933 \cite{thompson_likelihood_1933}, and despite subsequent research endeavors, it remains partially unresolved. Gittins made noteworthy advancements by addressing the class of index policies of the classical MAB problem, effectively reducing the complexity from an $N$-dimensional problem to~$N$ individual $1$-dimensional problems \cite{gittins_bandit_1979, gittins_dynamic_1979}.

Whittle further extended the classical MAB formulation and introduced the more comprehensive variant, the restless multi-armed bandit (RMAB) problem~\cite{whittle_restless_1988}. In RMAB, the player is granted to select~$M$ arms from the available $N$ arms (where $1\leq M\leq N$), and passive (unselected) arms can also alter their states, either of which generalizations made Gittins's approach (called Gittins index) suboptimal. Employing Lagrangian relaxation techniques, Whittle devised an indexing policy that generalizes Gittins index for a much broader spectrum of problems. This policy assigns an index (called Whittle index) to each arm dependent on the state of the arm. At each time, the player selects arms currently with the top~$M$ largest indices. Whittle's generalization has exhibited remarkable performance in both theoretical and numerical investigations \cite{brown_index_2020, chen_whittle_2021, weber_index_1990, zayas-caban_asymptotically_2019}. Nonetheless, establishing the essential condition for the existence of the Whittle index, known as indexability, and computing the Whittle index when it does exist pose significant challenges. Researchers have demonstrated that the RMAB problem with a finite state space is classified as PSPACE-hard problem \cite{papadimitriou_complexity_1994}. In this paper, we formulate the channel allocation problem as an RMAB with an infinite state space and construct an efficient algorithm to compute the Whittle index with arbitrary precision that achieves a near-optimal performance.

\subsection{Related Work}
Numerous prior studies have investigated the channel allocation problem in similar network models \cite{xie_spectrum_2012, lopez-perez_power_2014, anand_machine_2023, lin_interference_2020, eslami_joint_2023, famili_ofdra_2023, sharma_secrecy_2023, zhang_joint_2022, xu_robust_2022, li_distributed_2022, scopelliti_energy-quality_2020, pu_novel_2020}. Previous research has shown that the problem can be viewed as an MAB \cite{zhaokrish07ICC, zhaoEtal08TWC, ninomora09ICNCO, liu_distributed_2010, liuzhao10TSP, liu_indexability_2010, liuEtal10TSP, gai_distributed_2014, meshramEtal18TAC, kaza19TCCN, ninomora20MOR, liu_low-complexity_2022, ninomora26arXive}. Specifically, \cite{zhaokrish07ICC, zhaoEtal08TWC, ninomora09ICNCO, liu_indexability_2010, liuEtal10TSP, meshramEtal18TAC, kaza19TCCN, ninomora20MOR, liu_low-complexity_2022, ninomora26arXive} were all based on the two-state Gilbert-Elliott channel model shown in \figurename\ref{channel model} and formulated the RMAB as a partially observable Markov decision process (POMDP) with an infinite state space (belief space). In general, POMDP~\cite{sondik_optimal_1978} is numerically solved by dynamic programming and suffers from the curse of dimensionality. For the simplest case where channel state can be perfectly observed after sensing, Liu and Zhao~\cite{liu_indexability_2010} theoretically proved indexability and solved for the Whittle index in closed-form, leading to an efficient algorithm for heterogeneous channels (channels with different state transition probabilities and bandwidths) as a generalization of the myopic policy for homogeneous channels \cite{zhaokrish07ICC, zhaoEtal08TWC}. In \cite{ninomora09ICNCO,liuEtal10TSP, meshramEtal18TAC,kaza19TCCN, liu_low-complexity_2022}, the problem was considered with sensing errors and became fundamentally more complex than the perfect observation model considered in \cite{liu_indexability_2010}. This is because our belief space consists of the direct product of uncountable real intervals (each interval is exactly $[0,1]$ representing the conditional probability that the underlying channel state is~$1$ given the observation history) from all channels. The complex transition behavior in the belief space caused by sensing errors is the main difficulty in analyzing the dynamic programming equations. In~\cite{ninomora09ICNCO}, the class of threshold policies was assumed to simplify the indexability analysis and the numerical computation of the approximate Whittle index. Subsequent work \cite{meshramEtal18TAC, kaza19TCCN, liu_low-complexity_2022} gradually proved the optimality of threshold policies and solved for the approximate Whittle index in closed-form. For homogeneous channels (channels with same state transition probabilities and bandwidths), the Whittle index policy is reduced to the myopic policy with a simple structure and optimal performance under certain conditions~\cite{liuEtal10TSP}. 

Nevertheless, the imperfect observation model adopted in \cite{ninomora09ICNCO, liuEtal10TSP, meshramEtal18TAC,kaza19TCCN, liu_low-complexity_2022} was based on a binary feedback mechanism, e.g., ACK/NAK at the end of a transmission. In contrast, we consider the CQI feedback model as a generalization of the binary observation. The practical scenario where CQI is preferred to ACK/NAK can be found in Cellular Multicast/Broadcast Services (such as 5G MBS or LTE eMBMS) according to the modern wireless network standards \cite{3GPP_PHY,3GPP_UE}. From those standards, the main advantages of using CQI over ACKs are the following: 1. CQI provides finer information on channel qualities than ACK/NAK. Specifically, CQI informs the transmitter on the quality of the channel, allowing the change of power or modulation and coding scheme (MCS) to avoid band waste or save the power for data transmission if the channel condition is bad and move to a different band with better estimated CQI in the next time slot. ACK/NAK only tells the success or failure after a packet has been received and error checked with CRC (Cyclic Redundancy Check). The only thing that the transmitter can do is to retransmit again in HARQ (Hybrid Automatic Repeat Request); 2. Unlike ACK/NACK, which is an explicit receipt for a specific packet, CQI is decoupled from data packets and simply a measurement of the channel's signal-to-interference-plus-noise ratio (SINR). The user equipment (UE) measures this by listening to periodic reference signals (CSI-RS) broadcast by the tower. Here another practical example of using CQI over ACK/NACK is due to the round-trip propagation delay to a GEO satellite (roughly 500 to 600 milliseconds) \cite{3GPP_NTN}. Standard HARQ protocols require an ACK/NACK within a few milliseconds. By the time a NACK reaches the satellite, the timing window has expired, and holding millions of un-ACKed packets in memory for half a second exceeds the buffer capacity of the hardware. However, the receiver can estimate CQI even if it is in a connected state, by passively listening other signals sent by the transmitter to other receivers. It can send CQI feedback to the transmitter just as ACKs but with a wider range of information and without synchronization to transmission of a specific data packet; 3. ACK/NAK may be affected by temporary noise surge. CQI may be averaged over time and still provides stable channel information to the transmitter. When packet-by-packet feedback is impossible, the physical channel conditions (like "rain fade" or atmospheric attenuation) change relatively slowly. The ground terminal can reliably measure the channel and send CQI updates over a dedicated telemetry uplink. The satellite uses this CQI to perform adaptive coding and modulation (ACM), seamlessly lowering the data rate during a rainstorm to prevent packets from dropping in the first place. Given these advantages and practical significance of using CQI, we build our dynamic spectrum access model on top of the modern wireless system by allowing the sub-network base station to acquire the CSI feedback from the UE about the channel it currently chooses for potential data transmission. If the CSI report indicates a good channel condition (no traffic of the parent network or heavy noise), the station schedules the UE to transmit its own (sub-network) data. So the CQI report of a channel is requested by the station before the data transmission if the station plans to schedule this channel for data transmission in the current time slot (aperiodic mechanism specified in \cite{3GPP_UE}). Nevertheless, the station can also request CQI from any channel at any time, independent of any data transmission. A natural scenario is that the station asks for CQI periodically from a channel that has not been chosen for a long time to keep its channel information up to date \cite{3GPP_UE}. In this scenario, the station has more accurate and robust CQI reports (as they can be averaged within a short history) to make better decisions about which channel to select for data transmission. Therefore our model maps well to realistic cellular or cognitive-radio systems by using CQI reports for channel selections. However, due to the multiple levels in the CQI report, the belief transition behavior becomes more complex and theoretical analysis on the optimality of threshold policies, indexability and numerical computation of the Whittle index seemed impossible without very coarse approximations~\cite{elmaghraby_femtocell_2018}. 

In this paper, we overcome this challenge by establishing tight bounds on the value functions and their derivatives in dynamic programming under certain conditions, leading to the theoretical proof of the optimality of threshold policies, indexability and the closed-form Whittle index of geometrically decreasing approximation error with iterations. Furthermore, even if those conditions do not hold, our algorithm still produces an efficient index policy with strong performance demonstrated by numerical simulations. The methodology proposed in this paper also applies to other practical problems with multi-level feedback mechanisms. For example, in a cybersecurity monitoring system, the network controller selects a local component to monitor its traffic, query response, power usage data, etc. Different levels of information indicate the probabilities that the component is in the wrong state (malfunctional or under attack). The objective is to correctly target at those bad components so the repair team is not sent for nothing. For financial investments, different products are classified as `worth buying' or `not' in each decision epoch. The investor needs to estimate the probability that a product is in the `worth buying' state based on the market data as multi-level information (e.g., the stock price trend, the company's annual report, new technologies and competitions emerged).

Beyond DSA, the busy/idle channel model was also adopted in treatment adherence~\cite{ninomora26arXive}, throughput optimization for uncooperative users \cite{stahlbuhk20ToN}, and optimization of age-of-information problems for multiuser uplinks scheduling \cite{gongEtal20GLOBECOM, liuEtal23ToC}. In~\cite{ninomora26arXive}, a perfect observation model was assumed to satisfy the complex conditions proposed in~\cite{ninomora20MOR} on a verification theorem that tests the optimality of threshold policies and indexability simultaneously. In~\cite{stahlbuhk20ToN}, the central controller of the network needs to assign channels to adaptive (cooperative) users to transmit data. However, there are also uncooperative users who transmit data whenever their queues are nonempty and they do not inform the central controller their queue backlogs. So the central controller needs to estimate the queue backlogs of the uncooperative users to maximum the transmission opportunities for the adaptive users. The problem was thus formulated as a POMDP and a queue-length-based scheduling policy was analyzed to obtain the throughput-stability region of the network. In \cite{gongEtal20GLOBECOM, liuEtal23ToC}, the multiuser uplink system has one access point (AP) and a set of nodes. Each node will receive status update packets and maintain a local age of the last update packet. The goal of the AP is to select one node to transmit its update packet (synchronize its age with AP's record) such that the AP's record of weighted sum of ages of all nodes is minimized (ensuring information freshness). Because the AP does not know the update age of a node unless it transmits the information to the AP, the problem was formulated as a POMDP with a state space including the probability estimation of each node's local age and a low-complexity scheduling policy was proposed with strong performance. The case of observation errors was not considered in \cite{stahlbuhk20ToN, gongEtal20GLOBECOM, liuEtal23ToC}.

\subsection{Organization}
The rest of the paper is organized as follows. In Section~\ref{sec: system model and problem formulation}, we present the system model and the RMAB formulation of the channel allocation problem. In Section~\ref{sec: indexability and threshold policy}, we introduce the basic concepts of indexability and Whittle index, prove some important properties of value functions, and subsequently derive sufficient conditions for the threshold structure of the optimal policy and the indexability for the decoupled single-arm problem. In Section~\ref{sec: approximated whittle index}, we propose an iterative method to approximate value functions, provide 
the approximate Whittle index in closed-form, and finish the construction of the approximate Whittle index policy (AWI). In Section~\ref{sec: simulation results}, results of numerical experiments are provided to demonstrate the effectiveness of our generalized AWI compared with other heuristic policies. Finally, Section~\ref{sec: conclusion} concludes this paper with possible directions for future research.

\section{System Model and Problem Formulation}
\label{sec: system model and problem formulation}
\begin{figure}[!t]
\centering
\includegraphics[width=3.5in]{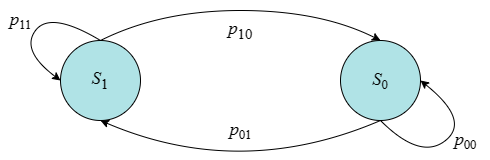}
\caption{The Gilbert-Elliott channel model}
\label{channel model}
\end{figure}

\subsection{System Model}
The system model comprises a parent network and a sub-network that shares a certain number of channels with the former. The goal of channel allocation is to optimize the selection of channels of the sub-network so that we can maximize total throughput while minimizing interference to the parent network.

The larger-sized MBS is a base station of the parent network while a smaller FBS of the sub-network is choosing among a set of channels allowed to be shared with the MBS. Due to resource limits and interference constraints, at each time, the FBS can only choose a portion of the shared channels to assign to the sub-network users for transmitting data. However, the data transmission can be successful only if the assigned channel is in the good state, currently not used by the MBS of the parent network. For each shared channel, we adopt the widely recognized two-state Markov model, also known as the Gilbert-Elliott model, as shown in \figurename\ref{channel model}. A channel can be in either the poor state or the good state.

The FBS does not directly observe the current channel states or CQI. Based on the CQI information previously reported by the sub-network users, the FBS chooses a portion of the shared channels and informs the users which channels are assigned for data transmission. By employing cognitive capabilities, a user can acquire the CQI from its assigned channel by listening to the channel without requiring knowledge of its actual state. The acquired CQI is then reported to the FBS and used as observation information for subsequent channel-selection decisions. Note that CQI levels are not channel states but contain information about the current state of a channel. Specifically, there is a known probability of observing a specific CQI level given a channel state.

\subsection{Restless Multi-armed Bandit Formulation}
For each shared channel, let $S \in \{S_0 := 0 \text{ (poor)}, S_1 := 1 \text{ (good)}\}$ represent the current state of a channel. The transition probabilities between the two states $0$ and $1$ are denoted by $p_{00}$, $p_{01}$, $p_{10}$ and $p_{11}$. The CQI acquired by the FBS from a sub-network user for a channel at time~$t$ is denoted by $q\left(t\right)\in\mathbb{Z}$. In our study, let~$K$ denote the number of CQI levels under investigation. Consequently, $q(t)$ satisfies the condition that $1\leq q\left(t\right)\leq K$. Our main notations are summarized in Table~\ref{tab:notations}.

Assume that the sub-network has~$N$ shared channels that are available to use. At each time instance, it is required to select~$M$ channels ($1\le M\le N$) for data transmission. For a given channel~$n$, if it is chosen (active) at time~$t$, we denote the action $a_n\left(t\right)=1$; otherwise $a_n\left(t\right)=0$. This is subject to the constraint that the number of active channels at each time equals~$M$ ($\sum_{n=1}^{N}{a_n\left(t\right)}=M$). At each decision epoch, the FBS selects the active channels using the current beliefs, which incorporate the CQI feedback from the preceding slot. An active channel is selected for transmission, whereas a passive channel is not selected and provides no CQI feedback.

Considering the unavailability of real-time channel state information in HetNets, we incorporate the general theory of partially observable Markov decision processes (POMDPs) \cite{sondik_optimal_1978} into our model. Specifically, we employ the belief state vector as the system state for decision-making purposes. The belief state vector, denoted by $\boldsymbol{\omega}(t)$, comprises the instantaneous belief states of the~$N$ channels. Specifically, the belief state of channel~$n$ at time~$t$ is defined as
\begin{equation}
    \omega_n(t) = \Pr(S_n(t) = 1 \mid \text{past observations on channel } n).
\end{equation}
Following the Bayes' rule, it can be proven that the evolution of belief state itself is a Markov process with an infinite state space. The belief is updated differently under the two actions. For an active channel, the CQI feedback is used to correct the belief by Bayes' rule, and the corrected belief is then propagated through the channel transition matrix. For a passive channel, no CQI is available, so the belief is propagated directly through the transition matrix. Thus, the transition probability of $\omega_n\left(t\right)$ can be described as
\begin{equation}
    \omega_n(t+1)=
    \begin{cases}
        \omega_{n,i}(t), & a_n(t)=1,\ q_n(t)=i \\
	\mathcal{T}_n(\omega_n(t)), & a_n(t)=0.
    \end{cases}
    \label{eq: 18}
\end{equation}
where $\mathcal{T}_n(\omega)$ is the belief update when passive action is taken on belief state $\omega$, hence no observation of CQI such that
\begin{equation}
    \mathcal{T}_n(\omega)=p^{(n)}_{11} \omega+p^{(n)}_{01}(1-\omega).
\end{equation}
We can also compute the belief update for successively~$k$ passive steps from the initial belief state~$\omega$ as
\begin{multline}
    \mathcal{T}_n^k(\omega)\\
    =\frac{p^{(n)}_{01}-\left(p^{(n)}_{11}-p^{(n)}_{01}\right)^k\left(p^{(n)}_{01}-\left(1+p^{(n)}_{01}-p^{(n)}_{11}\right) \omega\right)}{1+p^{(n)}_{01}-p^{(n)}_{11}}
\end{multline}
Moreover, if we let~$k$ tend to infinity, we can get the steady state belief value
\begin{equation}
    \omega_{n,s} = \lim_{k\to\infty}\mathcal{T}_n^k(\omega) = \frac{p^{(n)}_{01}}{1 + p^{(n)}_{01} - p^{(n)}_{11}}.
\end{equation}
Alternatively, when the channel is active and the observed CQI level is~$i$, the belief update is
\begin{equation}
\begin{aligned}
    \omega_{n}(t+1) &= \mathcal{T}(\Pr(S_n(t) = 1 \mid q_n(t) = i,\ \omega_n(t))) \\
    &= \mathcal{T}\left(\frac{p^{(n)}_{i,1}\omega_n(t)}{p^{(n)}_{i,1} \omega_n(t)+p^{(n)}_{i,0}(1-\omega_n(t))}\right) \\
    &= \frac{p^{(n)}_{11} p^{(n)}_{i,1} \omega_n(t)+p^{(n)}_{01} p^{(n)}_{i,0}(1-\omega_n(t))}{p^{(n)}_{i,1} \omega_n(t)+p^{(n)}_{i,0}(1-\omega_n(t))},
\end{aligned}
\end{equation}
where $p^{(n)}_{i,1}, p^{(n)}_{i,0}$ are the probabilities of observing $q_n(t) = i$ when the channel~$n$ is in the good ($S = 1$) or poor ($S = 0$) state, respectively. And it is obvious that $\sum_{i=1}^{K} p^{(n)}_{i,1}=1$ and $\sum_{i=1}^{K} p^{(n)}_{i,0}=1$. These probabilities characterize the CQI feedback generated by an active channel; no CQI observation is generated when the channel is passive. Moreover, the probability of observing a certain CQI given $\omega_n(t)$ is given by
\begin{equation}\label{eqn:fs_prob}
\begin{aligned}
    p_i^{(n)}(\omega_n(t)) &= \Pr(q_n(t) = i | \omega_n(t)) \\
    &= p_{i,1}^{(n)}\omega_n(t) + p_{i,0}^{(n)}(1 - \omega_n(t)).
\end{aligned}
\end{equation}

\begin{example}[3-state Feedback Model]
Here we illustrate a toy example to elaborate on the belief transition behaviors given the observation history. Let $K=3$ and
\[
\begin{pmatrix}
p^{(n)}_{1,1} & p^{(n)}_{2,1} & p^{(n)}_{3,1}\\
p^{(n)}_{1,0} & p^{(n)}_{2,0} & p^{(n)}_{3,0}\\
\end{pmatrix}
=\begin{pmatrix}
0.1 & 0.3 & 0.7 \\
0.7 & 0.3 & 0.1
\end{pmatrix},
\]
\[
\begin{pmatrix}
p^{(n)}_{00} & p^{(n)}_{01}\\
p^{(n)}_{10} & p^{(n)}_{11}
\end{pmatrix}
=\begin{pmatrix}
0.8 & 0.2\\
0.2 & 0.8
\end{pmatrix},
\]
\[
\omega_{n}(t)=0.5.
\]
By~\eqref{eqn:fs_prob} and Bayes' rule, if the observed CQI level is~$3$, then $\omega_{n}(t+1)=\mathcal{T}\left(\frac{0.5\cdot0.7}{0.5\cdot0.7+(1-0.5)\cdot0.1}\right)=\mathcal{T}(0.875)=0.875\cdot0.8+(1-0.875)\cdot0.2=0.725$. The updated values of the belief state given other observed CQI levels can be similarly computed.
\end{example}

Given the initial belief state vector $\boldsymbol{\omega}(1)$, we can formulate the channel allocation problem as a constrained optimization problem
\begin{align}
    \max_{\pi} &\mathbb{E}_\pi\left[\sum_{t=1}^\infty\beta^{t-1}R_\pi(t) \mid \boldsymbol{\omega}(1)\right], \\
    &\text{subject to } \sum_{n=1}^N a_n(t) = M,
\end{align}
where $\beta\in[0,1]$ is the discount factor to balance the importance of the instantaneous and future rewards. And the reward function is defined as
\begin{equation}
    R_\pi(t) = \sum_{n = 1}^N \mathbbm{1}_{\{a_n(t) = 1\}}S_n(t)B_n,
\end{equation}
where $B_n$ is the throughput of channel $n$.

Introducing the Lagrangian multiplier (subsidy for passivity)~$m$ and applying the Lagrangian relaxation as Whittle did in~\cite{whittle_restless_1988}, we can simplify the $N$-channel optimization problem to~a single-channel scenario:
\begin{equation}
    \max_{\pi: \omega_n(t) \rightarrow\{0,1\}} \mathbb{E}_\pi\left[\sum_{t=1}^{\infty} \beta^{t-1}\Tilde{R}_{\pi}^{(n)}(t) \mid \omega_n(1) = \omega\right],
    \label{eq: 20}
\end{equation}
where
\begin{equation}
    \Tilde{R}_{\pi}^{(n)}(t) = \mathbbm{1}_{\{a_n(t)=1\}} S_n(t) B_n+ m\cdot \mathbbm{1}_{\{a_n(t)=0\}}.
\end{equation}
The optimal value of the unconstrained optimization problem~(\ref{eq: 20}) is denoted by $V^{(n)}_{\beta,m}(\omega)$, and it is equivalent to
\begin{equation}
    V^{(n)}_{\beta, m}(\omega) = \max\{V^{(n)}_{\beta, m}(\omega; a = 0),\ V^{(n)}_{\beta, m}(\omega; a = 1)\},
\end{equation}
where $V^{(n)}_{\beta, m}(\omega; a = 1)$ and $V^{(n)}_{\beta, m}(\omega; a = 0)$ represent the optimal value of~(\ref{eq: 20}) whether channel~$n$ is chosen or not at the initial belief state~$\omega$, respectively. To simplify the presentation and without loss of generality, we will omit the superscript~$(n)$ and subscript~$n$ and set $B = 1$, considering a single-armed bandit problem in the following. And we can prove that the value function for passive and active actions satisfy the dynamic equations below
\begin{align}
    V_{\beta,m}(\omega; a=0) &= m + \beta V_{\beta,m}(\mathcal{T}(\omega)),\label{eq: 1} \\
    V_{\beta,m}(\omega; a=1) &= \omega + \beta\sum_{i=1}^{K}[p_{i,1} \omega+p_{i,0}(1-\omega)] V_{\beta,m}\left(\omega_i\right).\label{eq: 2}
\end{align}

For the case of perfect observation~\cite{liu_indexability_2010}, it is not hard to see that $V_{\beta,m}(\omega; a=1)$ is linear with~$\omega$. Together with the fact that $V_{\beta,m}(\omega; a=0)$ is convex in~$\omega$, the images of $V_{\beta,m}(\omega; a=1)$ and $V_{\beta,m}(\omega; a=0)$ have a unique intersecting point. In other words, the optimal policy for the single-armed bandit is a threshold policy, leading to closed-form solutions of the value functions and the Whittle index. However, in the presence of observation errors, both $V_{\beta,m}(\omega; a=1)$ and $V_{\beta,m}(\omega; a=0)$ are nonlinear and the analysis on the optimality of threshold policies and value functions becomes very difficult. To combat this challenge, we start to analyze the value functions over a finite time horizon where backward induction becomes possible. By carefully bounding the value functions and their derivatives along the backpropagation process, we are able to prove the optimality of threshold policies under certain conditions. Then these properties will be proved to still hold when the time horizon goes to infinity by the uniform convergence theorem. Based on the optimal threshold policy and the bounds on the value functions, we can further establish indexability by showing the monotonic property of the threshold with~$m$ and thus solve for the Whittle index of~$\omega$ defined as the minimum value of~$m$ that makes it the threshold, i.e., $V_{\beta,m}(\omega; a=1)=V_{\beta,m}(\omega; a=0)$. The detailed derivations are given in Sec.~\ref{sec: indexability and threshold policy} and Sec.~\ref{sec: approximated whittle index}.


\section{Indexability and Threshold Policy}
\label{sec: indexability and threshold policy}
Define passive set $P(m)$ as the collection of all belief states where the optimal action is to be passive ($a = 0$)
\begin{equation}
    P(m) = \{\omega:\ V_{\beta,m}(\omega; a=1) \leq V_{\beta,m}(\omega; a=0)\}.
\end{equation}
A restless multi-armed bandit is called {\it indexable} if for each single-armed bandit problem with Lagrangian multiplier~$m$, the passive set $P(m)$ monotonically expands from the empty set to the entire state space as~$m$ increases from $-\infty$ to $+\infty$~\cite{weber_index_1990}. According to~\cite{liu_indexability_2010}, under indexability, the Whittle index $W(\omega)$ for~a particular belief state $\omega$ is defined as follows:
\begin{equation}
    W(\omega) = \inf\{m:\ V_{\beta,m}(\omega; a=1) = V_{\beta,m}(\omega; a=0)\}.
\end{equation}

\subsection{Properties of Value Functions}
To prove the indexability of the RMAB problem and derive the threshold structure of the optimal policy for relaxed single-armed bandit problem~(\ref{eq: 20}), we need to investigate the properties of value functions $V_{\beta,m}(\omega; a=1)$, $V_{\beta,m}(\omega; a=0)$ and $V_{\beta,m}(\omega)$. We divide this process into two steps, following the approach in~\cite{liu_low-complexity_2022}. First, we examine the properties of value functions in finite horizons. Then, utilizing the uniform convergence theorem, we extend these conclusions to the case of an infinite horizon.

We introduce the $T$-horizon value function $V_{1,T,\beta,m}(\omega)$ as follows:
\begin{equation}
    V_{1,T,\beta,m}(\omega) = \max _{\pi} \mathbb{E}_\pi\left[\sum_{t=1}^{T} \beta^{t-1}\Tilde{R}_\pi(t) \mid \omega(1)=\omega\right],
\end{equation} where $\pi: \omega(t) \rightarrow\{0,1\}$ is a policy determining whether or not to activate the arm based on its current belief state.
Then it is obvious that
\begin{multline}
    V_{1,T,\beta,m}(\omega) \\
    = \max\{V_{1,T,\beta,m}(\omega; a = 0),\ V_{1,T,\beta,m}(\omega; a = 1)\}.
\end{multline}
And finite-horizon action value functions also satisfy the similar dynamic equations as (\ref{eq: 1}) and (\ref{eq: 2}):
\begin{align}
    & V_{1,T,\beta,m}(\omega; a = 0) = m + \beta V_{1,T-1,\beta,m}(\mathcal{T}(\omega)), \\
    & V_{1,T,\beta,m}(\omega; a = 1)\notag\\
    &\quad = \omega + \beta\sum_{i=1}^{K}[p_{i,1} \omega+p_{i,0}(1-\omega)] V_{1,T-1\beta,m}\left(\omega_i\right), \label{eq: 19}
\end{align}
where $V_{1,0,\beta,m}(\omega) \equiv 0$.
Utilizing the above recursive formulas and mathematical induction, we analyze the properties of $V_{1,T,\beta,m}(\omega)$.
\begin{lemma}
    $V_{1,T,\beta,m}(\omega)$ is piecewise linear and convex in both $\omega$ and $m$ for any $T \geq 1$.
    \label{lemma: 1}
\end{lemma}
\begin{IEEEproof}
    Consider first $T = 1$. It is clear that
    \begin{equation*}
        V_{1,1,\beta,m}(\omega) = \max\{m, \omega\} = \begin{cases}
            m, & \omega < m \\
            \omega, & \omega \geq m
        \end{cases}
    \end{equation*}
    is the maximum of two linear equations and thus piecewise linear and convex in both~$\omega$ and~$m$. Based on the recursive formulas and the induction hypothesis that $V_{1,T-1,\beta,m}(\omega)$ is piecewise linear in both~$\omega$ and~$m$, we can prove that $V_{1,T,\beta,m}(\omega;a=0)$ and $V_{1,T,\beta,m}(\omega;a=1)$ are still piecewise linear in both~$\omega$ and~$m$. Note that the recursive equation~(\ref{eq: 19}) leads to the following term in the expression of $V_{1,T,\beta,m}(\omega; a=1)$,
    \begin{equation*}
        \left[p_{i,1} \omega+p_{i,0}(1-\omega)\right] V_{1,T-1\beta,m}\left(\omega_i\right),
    \end{equation*}
    which has~a coefficient $\left[p_{i,1} \omega+p_{i,0}(1-\omega)\right]$ also appeared as the denominator in the expression of $\omega_i$. Finally we obtain that $V_{1,T,\beta,m}(\omega)$ is the maximum of two piecewise linear and convex functions and thus piecewise linear and convex in both~$\omega$ and~$m$.
\end{IEEEproof}

\begin{lemma}
    If $p_{11} > p_{01}$, $V_{1,T,\beta,m}(\omega)$ is monotonically increasing with $\omega\in [0,1]$ for any $T\geq 1$. 
    \label{lemma: 2}
\end{lemma}
\begin{IEEEproof}
    Since the function $V_{1,T,\beta,m}(\omega)$ is piecewise linear, demonstrating the monotonically increasing nature of the continuous function $V_{1,T,\beta,m}(\omega)$ with respect to $\omega$ can be achieved by proving
    \begin{equation}
        \label{eq: 3}
        V_{1,T,\beta,m}'(\omega) \geq 0, \quad \forall\ \omega\in (0,1),
    \end{equation}
    where $V_{1,T,\beta,m}'(\omega)$ represents the right derivative of the function $V_{1,T,\beta,m}(\omega)$ with respect to $\omega$. 
    
    In the case of $T = 1$, it is evident that $V_{1,1,\beta,m}(\omega) = \max\{\omega, m\}$ has a non-negative right derivative of either~$1$ or~$0$. Now, assuming that~(\ref{eq: 3}) holds for $T$, we analyse the case of $T + 1$. Let
    \begin{equation}
        V_{1,T+1,\beta,m}(\omega) = \max\{f_T(\omega), g_T(\omega)\}
    \end{equation}
    with
    \begin{align}
        f_T(\omega) &= m + \beta V_{1,T,\beta,m}(\mathcal{T}(\omega)), \\
        g_T(\omega) &= \omega + \beta\sum_{i=1}^K\left[ p_{i,1}\omega +  p_{i,0}(1 - \omega)\right]V_{1,T,\beta,m}\left(\omega_i\right).
    \end{align}
    Let 
    \begin{equation}
        \omega_i = \mathcal{T}\left(\frac{p_{i,1}\omega}{p_{i,1}\omega + p_{i,0}(1-\omega)}\right) = \mathcal{T}\left(\phi_i(\omega)\right).
    \end{equation}
    Taking the derivatives with respect to~$\omega$, we have
    \begin{align}
        f'_T(\omega) &= \beta(p_{11} - p_{01})V_{1,T,\beta, m}'(\mathcal{T}(\omega)),\label{eq: 4} \\	
        g'_T(\omega) &= 1 + \beta\sum_{i=1}^K(p_{i,1} - p_{i,0})V_{1,T,\beta,m}\left(\mathcal{T}\left(\phi_i(\omega)\right)\right)\notag \\
        & + \beta\sum_{i=1}^K V_{1,T,\beta,m}'\left(\mathcal{T}\left(\phi_i(\omega)\right)\right)\frac{p_{i,1}p_{i,0}(p_{11} - p_{01})}{p_{i,1}\omega + p_{i,0}(1-\omega)}.
        \label{eq: 5}
    \end{align}
    
    When $p_{11} > p_{01}$, according to the induction hypothesis and~(\ref{eq: 4}), it is straightforward to conclude that $f_T'(\omega) \geq 0$ and thus $f_T(\omega)$ is monotonically increasing. To prove the monotonicity of $g_T(\omega)$, we begin by analyzing the properties of the function $\phi_i$. Let~$P$ and~$N$ denote the sets of CQI signals that satisfy
    \begin{equation}
        P = \{i:\ p_{i,1} - p_{i,0}\geq 0\}, \quad N = \{i:\ p_{i,1} - p_{i,0} < 0\}.
            \label{eq:29}
    \end{equation}
    Suppose that $i\in P$ and $j\in N$, then we have
    \begin{multline}
        \phi_i(\omega) - \phi_j(\omega) \\
        = \frac{(p_{i,1}p_{j,0} - p_{i,0}p_{j,1})\omega(1-\omega)}{\left[p_{i,1}\omega + p_{i,0}(1-\omega)\right]\left[p_{j,1}\omega + p_{j,0}(1-\omega)\right]} \geq 0.
    \end{multline}
    Considering that~$\mathcal{T}$ and $V_{1,T,\beta,m}(\omega)$ are both monotonically increasing under the condition $p_{11} > p_{01}$, for any $i\in P$, $j\in N$, we obtain 
    \begin{equation}
        V_{1,T,\beta, m}\left(\mathcal{T}\left(\phi_i(\omega)\right)\right) \geq V_{1,T,\beta, m}\left(\mathcal{T}\left(\phi_j(\omega)\right)\right).
    \end{equation}
    Thus we have
    \begin{align*}
        g'_T(\omega) &\geq \beta\sum_{i\in P}(p_{i,1} - p_{i,0})V_{1,T,\beta,m}\left(\mathcal{T}\left(\phi_i(\omega)\right)\right) \\
        &+ \beta\sum_{j\in N}(p_{j,1} - p_{j,0})V_{1,T,\beta,m}\left(\mathcal{T}\left(\phi_j(\omega)\right)\right) \\
        &\geq \beta[\sum_{i\in P}(p_{i,1} - p_{i,0}) + \sum_{j\in N}(p_{j,1} - p_{j,0})] \\
        &~~~~\cdot\max_{j\in N}\{V_{1,T,\beta,m}\left(\mathcal{T}\left(\phi_j(\omega)\right)\right)\} \\
        &\geq 0.
    \end{align*}
    This establishes the property of $g_T(\omega)$ being monotonically increasing. Consequently, $V_{1,T,\beta,m}(\omega) = \max\{f_T(\omega), g_T(\omega\}$ is also monotonically increasing, thereby concluding the proof.
\end{IEEEproof}

\begin{assumption}
\label{ass:threshold}
For the relaxed single-armed bandit problem, the discount
factor $\beta\in(0,1)$ satisfies
\[
\beta \leq
\begin{cases}
\displaystyle
\frac{1}{
2(p_{11}-p_{01})
\left[
1+\sum_{i\in P}(p_{i,1}-p_{i,0})
\right]},
& p_{11}>p_{01},\\[3ex]
\displaystyle
\frac{1}{
(p_{01}-p_{11})
\left[
3+4\sum_{i\in P}(p_{i,1}-p_{i,0})
\right]},
& p_{01}>p_{11},
\end{cases}
\]
where $P$ is defined in~\eqref{eq:29}.
\end{assumption}

\begin{lemma}
    Under Assumption~\ref{ass:threshold}, the discount
factor also satisfies
    \begin{equation}
        \beta < \frac{1}{|p_{11} - p_{01}|\left[1 + 2\sum_{i\in P}(p_{i,1} - p_{i,0})\right]}.
    \end{equation}
    Then for all $T \geq 1$ and $\omega,\ \omega'\in[0,1]$, we have
    \begin{equation}
        |V_{1,T,\beta,m}(\omega) - V_{1,T,\beta,m}(\omega')| \leq C|\omega - \omega'|,
    \end{equation}
    where $C = \frac{1}{1 - \beta|p_{11} - p_{01}|\left[1 + 2\sum_{i\in P}(p_{i,1} - p_{i,0})\right]}$.
    \label{lemma: 3}
\end{lemma}
\begin{IEEEproof}
    In fact, we prove the conclusion above by directly proving that
    \begin{equation}
        |V_{1,T,\beta,m}'(\omega)| \leq C, \quad \forall\ T\geq 1,\ \omega\in(0,1).
    \end{equation}
    
    In the case of $T = 1$, $|V_{1,1,\beta,m}'(\omega)| \leq 1 \leq C$. Then under the induction hypothesis that $|V_{1,T,\beta,m}'(\omega)| \leq C$, we need to prove $|V_{1,T+1,\beta,m}'(\omega)| \leq C$. Recall the right derivatives $f_T'(\omega)$ and $g_T'(\omega)$ in~(\ref{eq: 4}) and~(\ref{eq: 5}). It is obvious that
    \begin{equation}
        \left|f'_T(\omega)\right| \leq \beta C |p_{11} - p_{01}|.
        \label{eq: 6}
    \end{equation}
    Meanwhile, note that
    \begin{equation}
        \left|\frac{p_{i,1}p_{i,0}(p_{11} - p_{01})}{p_{i,1}\omega + p_{i,0}(1-\omega)}\right| \leq \max\{p_{i,1}, p_{i,0}\}|p_{11} - p_{01}|,
    \end{equation}
    \begin{multline}
        \left|\sum_{i=1}^K\left(p_{i,1} - p_{i,0}\right)V_{1,T,\beta,m}\left(\mathcal{T}\left(\phi_i(\omega)\right)\right)\right| \\
        \leq C|p_{11} - p_{01}|\sum_{i\in P}(p_{i,1} - p_{i,0}).
    \end{multline}
    Thus we obtain the bound on $g'_T(\omega)$
    \begin{equation}
        |g'_T(\omega) - 1| \leq \beta C |p_{11} - p_{01}|\left[1 + 2\sum_{i\in P}(p_{i,1} - p_{i,0})\right].
        \label{eq: 7}
    \end{equation}
    From~(\ref{eq: 6}) and~(\ref{eq: 7}), we conclude that
    \begin{multline*}
        \left|V'_{1,T+1,\beta,m}(\omega)\right| \\
        \leq 1 + \beta C |p_{11} - p_{01}|\left[1 + 2\sum_{i\in P}(p_{i,1} - p_{i,0})\right] = C.
    \end{multline*}
    The induction process implies that
    \begin{equation*}
        |V_{1,T,\beta,m}'(\omega)| \leq C, \quad \forall\ T\geq 1,\ \omega\in(0,1).
    \end{equation*}
    Thus the proof is finished.
\end{IEEEproof}

\begin{lemma}
    Under Assumption~\ref{ass:threshold}, we have
\begin{equation}
        V'_{1,T,\beta,m}(\omega; a = 1) \geq V'_{1,T,\beta,m}(\omega; a = 0),
    \end{equation}
where $V'_{1,T,\beta,m}(\omega;a=k)$ denotes the right
derivative of $V_{1,T,\beta,m}(\omega;a=k)$ at $\omega$ for
$k\in\{0,1\}$.
    \label{lemma: 4}
\end{lemma}
\begin{IEEEproof}
    Again, we prove by mathematical induction on the time horizon~$T$. When $T = 1$, it is clear that
    \begin{equation*}
        V'_{1,1,\beta,m}(\omega; a=1) = 1 > V'_{1,1,\beta,m}(\omega; a=0) = 0.
    \end{equation*}
    Assume that $V'_{1,T,\beta,m}(\omega; a=1) \geq V'_{1,T,\beta,m}(\omega; a=0)$. If $p_{01} > p_{11}$ and
    \begin{equation*}
        \beta \leq \frac{1}{(p_{01} - p_{11})\left[3 + 4\sum_{i\in P}(p_{i,1} - p_{i,0})\right]},
    \end{equation*}
    we obtain that
    \begin{equation*}
        \beta C(p_{01} - p_{11}) \leq 1 - \beta C(p_{01} - p_{11})\left[1 + 2\sum_{i\in P}(p_{i,1} - p_{i,0})\right],
    \end{equation*}
    which shows that $f'_T(\omega) \leq g'_T(\omega)$ according to~(\ref{eq: 6}) and~(\ref{eq: 7}). 
    
    On the other hand, if $p_{11} > p_{01}$, $V_{1,T,\beta,m}(\omega)$ is increasing with $\omega$ with nonnegative right derivatives by Lemma \ref{lemma: 2}. We can thus obtain tighter bounds on $f'_T(\omega)$ and $g'_T(\omega)$
    \begin{equation*}
        1 \leq g'_T(\omega) \leq 1 + \beta C(p_{11} - p_{01})\left[1 + 2\sum_{i\in P}(p_{i,1} - p_{i,0})\right],
    \end{equation*}
    \begin{equation*}
        0 \leq f'_T(\omega) \leq \beta C(p_{11} - p_{01}).
    \end{equation*}
    When we choose
    \begin{equation*}
        \beta \leq \frac{1}{2(p_{11} - p_{01})\left[1 + \sum_{i\in P}(p_{i,1} - p_{i,0})\right]},
    \end{equation*}
    it is clear that $\beta C(p_{11} - p_{01}) \leq 1$, which shows that $f'_T(\omega) \leq g'_T(\omega)$. The proof is thus complete.	
\end{IEEEproof}

\subsection{Optimality of Threshold Policy}
In this section, we demonstrate that the optimal single-armed policy is~a threshold policy, subject to the constraints on the discount factor~$\beta$ outlined in the previous section. We still begin with the finite-horizon case.

For a $T$-horizon single-armed bandit problem, a threshold policy $\pi_T$ is defined by a time-dependent real number $\omega_{T,\beta}(m)$ such that
\begin{equation}
    a_{T,m}(\omega)= 
    \begin{cases}
        1, & \text { if } \omega>\omega_{T,\beta}(m);\\
        0, & \text { if } \omega \leq \omega_{T,\beta}(m).
    \end{cases}
\end{equation}
Intuitively, as the value of~$\omega$ increases, so does the expected immediate reward, which in turn makes activating the arm more appealing. The following theorem formalizes this intuition under specific conditions.
\begin{theorem}
    Under Assumption~\ref{ass:threshold}, for any $T \geq 1$, the optimal $T$-horizon single-armed policy $\pi^*_T$ is a threshold policy, which means that there exists $\omega^*_{T,\beta}(m)\in\mathbb{R}$ such that under $\pi^*_T$, the optimal action is 
    \begin{equation}
        a^*_{T,m}(\omega)= 
        \begin{cases}
            1, & \text { if } \omega>\omega^*_{T,\beta}(m);\\
            0, & \text { if } \omega \leq \omega^*_{T,\beta}(m).
        \end{cases}
    \end{equation}
    Furthermore, at the threshold belief state $\omega^*_{T,\beta}(m)$,
    \begin{equation}
        V_{1,T,\beta,m}(\omega^*_{T,\beta}(m); a=1) = V_{1,T,\beta,m}(\omega^*_{T,\beta}(m); a=0).
    \end{equation}
    
    \label{theorem: 1}
\end{theorem}

In the next theorem, we show that the optimal single-armed policy over the infinite horizon is also a threshold policy under the same conditions.
\begin{theorem}
    Under Assumption~\ref{ass:threshold}, fix the Lagrangian multiplier~$m$. The finite-horizon value functions $V_{1,T,\beta,m}(\cdot)$, $V_{1,T,\beta,m}(\cdot; a=1)$ and $V_{1,T,\beta,m}(\cdot; a=0)$ uniformly converge to the infinite-horizon value functions $V_{\beta,m}(\cdot)$, $V_{\beta,m}(\cdot; a=1)$ and $V_{\beta,m}(\cdot; a=0)$ which consequently possess the same properties established in Lemma \ref{lemma: 2}-\ref{lemma: 3} and Theorem~\ref{theorem: 1}.
    \label{theorem: 2}
\end{theorem}

Utilizing the properties of value functions, the proof of Theorem~\ref{theorem: 1} and~\ref{theorem: 2} are the same as those of Theorem~2.5 and~2.6 in~\cite{liu_low-complexity_2022}, so we omit it here. Thus far we have established the threshold structure of the optimal single-armed policy based on the analysis of $V_{\beta,m}(\omega)$ as a function of the belief state~$\omega$ with~$m$ fixed.

\subsection{Indexability}
Based on the definition of indexability and the threshold structure of the optimal policy for the relaxed single-armed bandit problem~(\ref{eq: 20}), establishing the indexability of our model is equivalent to proving that the threshold $\omega^*_\beta(m)$ is monotonically increasing with respect to~$m$.

To investigate the sufficient condition for indexability, we now examine the properties of $V_{\beta,m}(\omega)$ as a function of the Lagrangian multiplier~$m$ with the initial belief state~$\omega$ fixed.

\begin{lemma}
    Given the initial belief state~$\omega$, value function $V_{\beta,m}(\omega)$ is convex in~$m$. Furthermore, the left and right derivatives of $V_{\beta,m}(\omega)$ with respect to~$m$ exist at every point $m_0\in\mathbb{R}$.
    \label{lemma: 5}
\end{lemma}
\begin{IEEEproof}
    In fact, for any given values of $m_1$, $m_2$ and $\theta\in(0,1)$, when we apply the optimal policy $\pi^*_\beta(\theta m_1 + (1-\theta)m_2)$ of $V_{\beta,\theta m_1 + (1-\theta)m_2}(\omega)$ to the problem with Lagrangian multipliers $m_1$ and $m_2$ respectively, it cannot surpass the performance achieved by the optimal policies $\pi^*_{\beta}(m_1)$ and $\pi^*_{\beta}(m_2)$ for $V_{\beta,m_1}(\omega)$ and $V_{\beta,m_2}(\omega)$. To be more explicit, let $r_a$ and $r_p(m)$ represent the expected total discounted reward from active and passive actions under policy $\pi^*_\beta(\theta m_1 + (1-\theta)m_2)$ applied to the problem with Lagrangian multipliers~$m$ and initial belief state~$\omega$, then we have
\begin{align*}
    \theta V_{\beta,m_1}(\omega) +& (1-\theta)V_{\beta,m_2}(\omega) \\
    \geq& \theta(r_a + r_p(m_1)) + (1-\theta)(r_a + r_p(m_2)) \\
    =& r_a + r_p(\theta m_1 + (1-\theta)m_2) \\
    =& V_{\beta,\theta m_1 + (1-\theta)m_2}(\omega).
\end{align*}
This shows the convexity of value function $V_{\beta,m}(\omega)$ with respect to~$m$. Since $V_{\beta,m}(\omega)$ is convex in~$m$, it is obvious that its left and right derivatives with~$m$ exist at every point $m_0\in\mathbb{R}$, according to the properties of convex functions.
\end{IEEEproof}

\begin{lemma}
    Given the initial belief state~$\omega$, value function $V_{\beta,m}(\omega)$ is differentiable almost everywhere in~$m$.
\end{lemma}
\begin{IEEEproof}
    Consider two policies $\pi^*_\beta(m_1)$ and $\pi^*_\beta(m_2)$ achieving $V_{\beta,m_1}(\omega)$ and $V_{\beta,m_2}(\omega)$ for any~$m_1$,~$m_2\in\mathbb{R}$, respectively. Utilizing the similar trick as the proof of Lemma~\ref{lemma: 5}, let~$r_a$ be the expected total discounted reward from the active action and $r_p(m)$ that from the passive action under $\pi^*_\beta(m_1)$ applied to the problem with Lagrangian multiplier~$m$, then
\begin{equation*}
    V_{\beta,m_1}(\omega) = r_a + r_p(m_1),\quad V_{\beta,m_2}(\omega) \geq r_a + r_p(m_2).
\end{equation*}	
Thus we can obtain that
\begin{align*}
    V_{\beta,m_1}(\omega) - V_{\beta,m_2}(\omega) &\leq r_a + r_p(m_1) - r_a - r_p(m_2) \\
    & = r_p(m_1 - m_2) \\
    & \leq \frac{1}{1-\beta}|m_1 - m_2|.
\end{align*}
Interchanging $V_{\beta,m_1}(\omega)$ and $V_{\beta,m_2}(\omega)$ by the symmetry of the problem, we obtain the conclusion that $V_{\beta,m}(\omega)$ is Lipschitz continuous in~$m$, that is
\begin{equation*}
    \left|V_{\beta,m_1}(\omega) - V_{\beta,m_2}(\omega)\right| \leq \frac{1}{1-\beta}|m_1 - m_2|.
\end{equation*}
Consequently, according to Rademacher theorem \cite{zajicek_elementary_1992}, $V_{\beta,m}(\omega)$ is differentiable almost everywhere in~$m$.
\end{IEEEproof}

In the following theorem, we formalize the relationship between value function and passive time and provide a sufficient condition for the indexability of our model.
\begin{theorem}
    Let $\Pi^*_\beta(m)$ denote the set of all optimal single-armed policies achieving $V_{\beta,m}(\omega)$ with initial belief state~$\omega$. Define the passive time as
    {\small\begin{equation}
        D_{\beta, m}(\omega) = \max _{\pi_\beta^*(m) \in \Pi_\beta^*(m)} \mathbb{E}_{\pi_\beta^*(m)}\left[\sum_{t=1}^{\infty} \beta^{t-1} \mathbbm{1}_{\{a(t)=0\}} \mid \omega(1)=\omega\right].
    \end{equation}}
    The right derivative of the value function $V_{\beta,m}(\omega)$ with~$m$, denoted by $\frac{d V_{\beta, m}(\omega)}{(d m)^{+}}$, exists at every value of~$m$ and
    \begin{equation}
        \left.\frac{d V_{\beta, m}(\omega)}{(d m)^{+}}\right|_{m=m_0} = D_{\beta, m_0}(\omega).
    \end{equation}
    Furthermore, the single-armed bandit is indexable if at least one of the following condition is satisfied:
    \begin{enumerate}
        \item for any $m_0\in[0,1)$, the optimal policy is~a threshold policy with threshold $\omega^*_\beta(m_0)\in[0,1)$ (if the threshold is~a closed interval then the right end is selected) and
        \begin{multline}
            \left.\frac{d V_{\beta, m}\left(\omega_\beta^*\left(m_0\right) ; a=0\right)}{(d m)^{+}}\right|_{m=m_0}\\
            >\left.\frac{d V_{\beta, m}\left(\omega_\beta^*\left(m_0\right) ; a=1\right)}{(d m)^{+}}\right|_{m=m_0}.
        \end{multline}
        \item for any $m_0\in\mathbb{R}$ and $\omega\in P(m_0)$, we have
        \begin{equation}
            \left.\frac{d V_{\beta, m}\left(\omega; a=0\right)}{(d m)^{+}}\right|_{m=m_0} \geq \left.\frac{d V_{\beta, m}\left(\omega; a=1\right)}{(d m)^{+}}\right|_{m=m_0}.
            \label{eq: 10}
        \end{equation}
    \end{enumerate}
    \label{theorem: 3}
\end{theorem}
The proof follows similarly from the argument of Theorem~1 in~\cite{liu_rmab_2025}. So we omit it here. 

Theorem~\ref{theorem: 3} provides a direct way for checking the indexability of the bandit problem with the help of the passive times. And because of this, we can present the sufficient condition for the indexability of our problem.

\begin{assumption}
\label{ass:indexability}
The discount factor satisfies $0<\beta\leq\frac{1}{2}$.
\end{assumption}

\begin{corollary}
    Under Assumption~\ref{ass:threshold} and Assumption~\ref{ass:indexability}, the restless single-armed bandit problem is indexable.
    \label{corollary: 1}
\end{corollary}
\begin{IEEEproof}
    Similar to the dynamic equations that value functions $V_{\beta, m}(\omega)$ satisfy, the passive time also has its own dynamic equations
    \begin{align}
        D_{\beta,m}(\omega;a=0) &= 1 + \beta D_{\beta,m}(\mathcal{T}(\omega)), \\
        D_{\beta,m}(\omega;a=1) &= \beta\sum_{i=1}^K[p_{i,1}\omega + p_{i,0}(1-\omega)]D_{\beta,m}(\omega_i).
    \end{align}
    Thus the equation~(\ref{eq: 10}) is equivalent to
    \begin{equation}
        1 + \beta D_{\beta,m}(\mathcal{T}(\omega)) \geq \beta\sum_{i=1}^K[p_{i,1}\omega + p_{i,0}(1-\omega)]D_{\beta,m}(\omega_i).
    \end{equation}
    The above inequality clearly holds if $\beta \leq 0.5$ since $D_{\beta,m}(\cdot)\in[0,\frac{1}{1-\beta}]$ for any $m\in\mathbb{R}$.
\end{IEEEproof}
 
\section{Approximated Whittle Index}
\label{sec: approximated whittle index}
According to Theorem~\ref{theorem: 2} and Corollary~\ref{corollary: 1}, we assume that the following condition is satisfied such that the optimal policy for the relaxed single-armed bandit problem is~a threshold policy and the indexability holds
\begin{equation}
    \label{eq: 15}
    \beta \leq \begin{cases}
        \min\left\{\frac{1}{2|p_d|\left[1 + \sum_{i\in P}(p_{i,1} - p_{i,0})\right]}, 0.5\right\}, & p_{11} > p_{01} \\
        \min\left\{\frac{1}{|p_d|\left[3 + 4\sum_{i\in P}(p_{i,1} - p_{i,0})\right]}, 0.5\right\}, & p_{11} < p_{01} \\
    \end{cases},
\end{equation}
where $p_d = p_{11} - p_{01}$. Given any belief state~$\omega$, to solve for the Whittle index $W(\omega)$ under indexability, we need to find out the minimum Lagrangian multiplier~$m$ that satisfies the system of equations below
\begin{equation}
    \left\{\begin{aligned}
        & V_{\beta,m}(\omega; a = 1) = V_{\beta,m}(\omega; a = 0), \\
        & V_{\beta,m}(\omega; a=1) = \omega + \beta\sum_{i=1}^{K}p_i(\omega) V_{\beta,m}\left(\omega_i\right), \\
        & V_{\beta,m}(\omega; a=0) = m + \beta V_{\beta,m}(\mathcal{T}(\omega)).
    \end{aligned}\right.
    \label{eq: 14}
\end{equation}

Before solving the equations above and deriving the approximated Whittle index utilizing the threshold structure of the optimal policy, we first present the concept of first crossing time. Given two belief states~$\omega$ and~$\omega'$, the first crossing time is defined as
\begin{equation}
    L(\omega, \omega') = \min_{0\leq k < \infty} \{k: \mathcal{T}^k(\omega) > \omega'\},
    \label{eq: 11}
\end{equation}
where we set $\mathcal{T}^0(\omega) = \omega$ and
\begin{equation*}
    L(\omega, \omega') = +\infty,\quad \text{if $\mathcal{T}^k(\omega) \leq \omega'$ for all $k \geq 0$.}
\end{equation*}
It is evident that $L(\omega, \omega')$ is the minimum time slots required for a belief state~$\omega$ to remain in the passive set $P(m)$ before the channel is chosen, given a threshold $\omega'\in[0,1)$. 

Consider the nontrivial case $p_{01}$, $p_{11}\in(0,1)$ and $p_{01}\neq p_{11}$ where the Markov chain of the internal arm states \{S(t)\} is aperiodic and irreducible and that the belief update is action-dependent. According to~(\ref{eq: 11}), we can figure out that, if $p_{11} > p_{01}$,
\begin{multline}
    L\left(\omega, \omega^{\prime}\right)\\
    = \begin{cases}
        0, & \omega>\omega^{\prime} \\
        \left\lfloor\log _{p_d}\frac{p_{01}-\omega^{\prime}\left(1-p_d\right)}{p_{01}-\omega\left(1-p_d\right)}\right\rfloor+1, & \omega \leq \omega^{\prime}<\omega_s \\
        \infty, & \omega \leq \omega^{\prime}, \omega^{\prime} \geq \omega_s
    \end{cases},
\end{multline}
where steady belief state $\omega_s = p_{01} / (1 + p_{01} - p_{11})$, else if $p_{11} < p_{01}$,
\begin{equation}
    L\left(\omega, \omega^{\prime}\right)= \begin{cases}
        0, & \omega>\omega^{\prime} \\ 
        1, & \omega \leq \omega^{\prime} < \mathcal{T}(\omega) \\ 
        \infty, & \omega\leq\omega^{\prime}, \mathcal{T}(\omega) \leq\omega^{\prime}
    \end{cases}.
\end{equation}
Using the first crossing time, the value function $V_{\beta,m}(\omega)$ can be expanded as
\begin{equation}
    V_{\beta,m}(\omega) = b_1(\omega) m + b_2(\omega) \Omega(\omega) \\
    + \sum_{i=1}^K b_{3,i}(\omega)V_{\beta,m}(f_i(\omega)),
    \label{eq: 12}
\end{equation}
where
\begin{align}
    & \Omega(\omega) = \mathcal{T}^{L(\omega, \omega^*_\beta(m))}(\omega), \\
    & b_1(\omega) = \frac{1-\beta^{L\left(\omega, \omega_\beta^*(m)\right)}}{1-\beta}, \\
    & b_2(\omega) = \beta^{L\left(\omega, \omega_\beta^*(m)\right)}, \\
    & b_{3,i}(\omega) = \beta^{L\left(\omega, \omega_\beta^*(m)\right) + 1}p_i\left(\mathcal{T}^{L\left(\omega, \omega_\beta^*(m)\right)}(\omega)\right),
\end{align}
and we define the function $f_i(\omega)$ as
\begin{equation}
    f_i(\omega) = \mathcal{T}\left(\phi_i\left(\mathcal{T}^{L\left(\omega, \omega_\beta^*(m)\right)}(\omega)\right)\right),
\end{equation}
while we take the notation that
\begin{equation}
    f_{i_2,i_1}(\omega) = f_{i_2} \circ f_{i_1}(\omega) = f_{i_2}(f_{i_1}(\omega)),\ i_1,i_2\in\{1,2,\cdots,K\}.
\end{equation}
The summation term within (\ref{eq: 12}) poses difficulties in solving for $V_{\beta,m}(\omega)$ since new belief states are introduced as unknowns. To tackle this difficulty, we approximate $V_{\beta,m}(\omega)$ in an iterative fashion. We compute the expanded form of $V_{\beta,m}(\omega)$, $V_{\beta,m}(f_{i_1}(\omega))$, $V_{\beta,m}(f_{i_2,i_1}(\omega))$, $\cdots$, $V_{\beta,m}(f_{i_n,\cdots,i_2,i_1}(\omega))$ one by one. In this way, we get the following sequence of equations
\begin{equation}
\begin{aligned}
    & V_{\beta, m}(\omega) = b_1(\omega) m + b_2(\omega)\Omega(\omega) + \sum_{i_1=1}^K b_{3,i_1}(\omega)V_{\beta,m}(f_{i_1}(\omega)), \\
    & V_{\beta, m}(f_{i_1}(\omega)) = b_1(f_{i_1}(\omega)) m + b_2(f_{i_1}(\omega))\Omega(f_{i_1}(\omega)) \\
    &\quad + \sum_{i_2=1}^Kb_{3,i_2}(f_{i_1}(\omega))V_{\beta,m}(f_{i_2,i_1}(\omega)),\, i_1 \in\{1,2,\cdots,K\} \\
    & V_{\beta, m}(f_{i_2, i_1}(\omega)) =  b_1(f_{i_2, i_1}(\omega)) m \\
    &\quad + b_2(f_{i_2, i_1}(\omega))\Omega(f_{i_2, i_1}(\omega)) \\
    &\quad + \sum_{i_3=1}^Kb_{3,i_3}(f_{i_2, i_1}(\omega))V_{\beta,m}(f_{i_3,i_2,i_1}(\omega)),\\
    &\quad\quad i_1, i_2 \in\{1,2,\cdots,K\} \\
    &\quad \cdots \\
    & V_{\beta, m}(f_{i_n,\cdots,i_1}(\omega)) = b_1(f_{i_n,\cdots,i_1}(\omega)) m \\
    &\quad + b_2(f_{i_n,\cdots,i_1}(\omega))\Omega(f_{i_n,\cdots,i_1}(\omega)) \\
    &\quad + \sum_{i_{n+1}=1}^Kb_{3,i_{n+1}}(f_{i_n,\cdots,i_1}(\omega))V_{\beta,m}(f_{i_{n+1},\cdots,i_1}(\omega)),\\
    &\quad\quad i_1,\cdots,i_n \in\{1,2,\cdots,K\}
\end{aligned}
\label{eq: 13}
\end{equation}
For sufficiently large iterative steps~$n$, we can get an estimation of $V_{\beta,m}(\omega)$ with an arbitrarily small error by setting
\begin{equation*}
V_{\beta,m}\left(f_{i_{n+1},\cdots,i_1}(\omega)\right) = 0,\ \forall\ i_1,\cdots,i_{n+1} = 1,2,\cdots,K.
\end{equation*}
During the computation of $V_{\beta, m}(f_{i_n,\cdots,i_2, i_1}(\omega))$, the error of this estimate is discounted by the factor $\beta$. As~a result, the backward computation process for $V_{\beta,m}(\omega)$ experiences a geometrically decreasing error propagation. Consequently, we can obtain an approximation of $V_{\beta,m}(\omega)$ with arbitrary precision for any $\omega\in[0,1]$, denoted as $\widehat{V}_{\beta,m,n}(\omega)$, where~$n$ represents the iteration steps here. See Theorem~\ref{thm:complexity} for details.

In conclusion, the $n$-iteration Whittle index is based on the solution of the system of equations (\ref{eq: 14}). To be more explicit, for any channel in belief state $\omega$, substituting~$\omega_i$ and $\mathcal{T}(\omega)$ respectively for~$\omega$ in the above system of equations (\ref{eq: 13}), we obtain $n$-iteration estimates of $V_{\beta,m}(\omega_i)$ and $V_{\beta,m}(\mathcal{T}(\omega))$. Thus we can use them in system of equations (\ref{eq: 14}) to compute the approximated Whittle index $W(\omega)$ by setting $\omega^*_\beta(m) = \omega$ according to the first equation of the system and Theorem \ref{theorem: 1}.

Moreover, in the case of large values of $\beta\in(0,1)$, where one or both of Assumptions~\ref{ass:threshold} and \ref{ass:indexability} may be violated (i.e., condition (\ref{eq: 15}) is not satisfied), we can still utilize the aforementioned process to find out the Lagrangian multiplier~$m$ that satisfies~(\ref{eq: 14}), if such~a solution exists. It is important to note that the approximated value functions $\widehat{V}_{\beta,m}(\omega;a=1)$ and $\widehat{V}_{\beta,m}(\omega;a=0)$ are linear in~$m$. The equality of these approximated value functions provides~a unique solution for~$m$, if it exists. This obtained~$m$, if it is indeed~a solution, can then be employed as the approximated Whittle index $W(\omega)$, without necessitating indexability or the threshold structure of the optimal policy. On the other hand, if it does not exist, we can simply set
\begin{equation}
    W(\omega) = \omega B.
\end{equation}

Before computing the closed form of $n$-iteration approximated Whittle index, we solve for the simplest case of $0$-iteration, which is referred to as the imperfect Whittle index. Setting $V_{\beta,m}(f_{i_1}(\omega)) = 0$, we can directly solve for $V_{\beta,m}(\omega)$ in closed form and thus obtain the estimates of $V_{\beta,m}(\omega_i)$ and $V_{\beta,m}(\mathcal{T}(\omega))$
\begin{align*}
    \widehat{V}_{\beta,m,0}(\omega_i) &= b_1(\omega_i) m + b_2(\omega_i) \Omega(\omega_i), \\
    \widehat{V}_{\beta,m,0}(\mathcal{T}(\omega)) &= b_1(\mathcal{T}(\omega)) m + b_2(\mathcal{T}(\omega)) \Omega(\mathcal{T}(\omega)).
\end{align*}
Then plugging them into the equations (\ref{eq: 14}), we can get~a simple linear equation with respect to~$m$:
\begin{equation*}
    c_1 m = c_0,
\end{equation*}
where
\begin{align}
    c_0 &= \omega + \beta\sum_{i=1}^Kp_i(\omega)\left[b_2(\omega_i)\Omega(\omega_i) - b_2(\mathcal{T}(\omega))\Omega(\mathcal{T}(\omega))\right], \\
    c_1 &= 1 + \beta\sum_{i=1}^Kp_i(\omega)\left[b_1(\mathcal{T}(\omega)) - b_1(\omega_i)\right].
\end{align}
If $c_1 \neq 0$, we can obtain the imperfect Whittle index as below
\begin{equation}
    \widehat{W}_0(\omega) = \left.\frac{c_0}{c_1}\right|_{\omega^*_\beta(m) = \omega}.
    \label{eq: 17}
\end{equation}

\subsection{Closed Form of Approximated Whittle Index}
\begin{theorem}
    Given iteration step~$n$, belief state $\omega$ and Lagrangian multiplier~$m$, setting $V_{\beta,m}(f_{i_{n+1},\ldots,i_1}(\omega)) = 0$ for any $i_1,\cdots,i_{n+1}\in\{1,2,\cdots,K\}$ in equation set~(\ref{eq: 13}), we get the $n$-iteration estimate of $V_{\beta,m}(\omega)$, $\widehat{V}_{\beta,m,n}(\omega)$, for $n \geq 0$ as below
    \begin{equation}
        \widehat{V}_{\beta,m,n}(\omega) = k_n(\omega)m + a_n(\omega),
        \label{eq: 16}
    \end{equation}
    where 
    \begin{equation}
        \begin{aligned}
            & k_0(\omega) = b_1(\omega), \\
            & k_n(\omega) = k_0(\omega) + \sum_{i_1=1}^K b_{3,i_1}(\omega)k_{n-1}(f_{i_1}(\omega)),\, n \geq 1
        \end{aligned}
    \end{equation}
    and 
    \begin{equation}
        \begin{aligned}
            & a_0(\omega) = b_2(\omega)\Omega(\omega), \\
            & a_n(\omega) = a_0(\omega) + \sum_{i_1=1}^K b_{3,i_1}(\omega)a_{n-1}(f_{i_1}(\omega)),\, n \geq 1.
        \end{aligned}
    \end{equation}
    Thus approximating value functions in this way, solving the system of equations (\ref{eq: 14}) and letting $\omega^*_\beta(m) = \omega$, we get the $n$-iteration approximated Whittle index, $\widehat{W}_n(\omega)$, as below if $1 + \beta\left(k_{n,0} - \sum_{i=1}^Kp_i(\omega)k_{n,i}\right) \neq 0$:
    \begin{equation}
        \widehat{W}_n(\omega) = \left.\frac{\omega + \beta\left(\sum_{i=1}^Kp_i(\omega)a_{n,i} - a_{n,0}\right)}{1 + \beta\left(k_{n,0} - \sum_{i=1}^Kp_i(\omega)k_{n,i}\right)}\right|_{\omega^*_\beta(m) = \omega}.
    \end{equation}
    where $k_{n,0} = k_n(\mathcal{T}(\omega))$, $a_{n,0} = a_n(\mathcal{T}(\omega))$ and $k_{n,i} = k_n(\omega_i)$, $a_{n,i} = a_n(\omega_i)$ for $i = 1,2,\cdots, K$.
\end{theorem}
\begin{IEEEproof}
    We first prove the equation~(\ref{eq: 16}) by mathematical induction. We start from the $0$-iteration case. Let $V_{\beta,m}(f_{i_1}(\omega)) = 0$ for any $i_1\in\{1,2,\cdots,K\}$ in the system of equations~(\ref{eq: 13}), thus it is easy to compute that
    \begin{equation*}
        \widehat{V}_{\beta,m,0}(\omega) = k_0(\omega)m + a_0(\omega)
    \end{equation*}
    where
    \begin{equation*}
        \begin{aligned}
            & k_0(\omega) = \frac{1 - \beta^{L(\omega,\omega^*_\beta(m))}}{1 - \beta} = b_1(\omega), \\
            & a_0(\omega) = \beta^{L(\omega,\omega^*_\beta(m))}\mathcal{T}^{L(\omega,\omega^*_\beta(m))}(\omega) = b_2(\omega)\Omega(\omega).
        \end{aligned}
    \end{equation*}
    Then in the $1$-iteration case, we set $V_{\beta,m}\!\left(f_{i_2,i_1}(\omega)\right)=0$ for any $i_1,i_2\in\{1,2,\cdots,K\}$. Solving the equation set~(\ref{eq: 13}), we get
    \begin{equation*}
        \widehat{V}_{\beta,m,1}(\omega) = k_1(\omega)m + a_1(\omega)
    \end{equation*}
    where
    \begin{align*}
        k_1(\omega) &= b_1(\omega) + \sum_{i_1 = 1}^K b_{3,i_1}(\omega)b_1(f_{i_1}(\omega)) \\
        &= k_0(\omega) + \sum_{i_1 = 1}^K b_{3,i_1}(\omega)k_0(f_{i_1}(\omega)), \\
        a_1(\omega) &= b_2(\omega)\Omega(\omega) + \sum_{i_1 = 1}^K b_{3,i_1}(\omega)b_2(f_{i_1}(\omega))\Omega(f_{i_1}(\omega)) \\
        &= a_0(\omega) + \sum_{i_1 = 1}^K b_{3,i_1}(\omega)a_0(f_{i_1}(\omega)).
    \end{align*}
    
    By mathematical induction, we assume that $\widehat{V}_{\beta,m,n}(\omega)$ satisfies the conclusion for any belief state $\omega$. Based on this, we compute the $(n+1)$-iteration estimate of $V_{\beta,m}(\omega)$. Note that $\widehat{V}_{\beta,m,n+1}(f_{i_1}(\omega))$ is the $n$-iteration estimate of $V_{\beta,m}(f_{i_1}(\omega))$ and is equal to 
    \begin{equation*}
        \widehat{V}_{\beta,m,n+1}(f_{i_1}(\omega)) = k_n(f_{i_1}(\omega))m + a_n(f_{i_1}(\omega)).
    \end{equation*}
    Then we have that
    \begin{align*}
        &\widehat{V}_{\beta,m,n+1}(\omega) \\
        =& b_1(\omega)m + b_2(\omega)\Omega(\omega) + \sum_{i_1=1}^Kb_{3,i_1}(\omega)\widehat{V}_{\beta,m,n+1}(f_{i_1}(\omega)) \\
        =& b_1(\omega)m + b_2(\omega)\Omega(\omega) \\
        &+ \sum_{i_1=1}^Kb_{3,i_1}(\omega)\left[k_n(f_{i_1}(\omega))m + a_n(f_{i_1}(\omega))\right] \\
        =& \left[b_1(\omega) + \sum_{i_1 = 1}^K b_{3,i_1}(\omega)k_n(f_{i_1}(\omega))\right]m \\
        &+ b_2(\omega)\Omega(\omega) + \sum_{i_1=1}^Kb_{3,i_1}(\omega)a_n(f_{i_1}(\omega)) \\
        =& k_{n+1}(\omega) m + a_{n+1}(\omega).
    \end{align*}
    Thus the proof of equation~(\ref{eq: 16}) is finished.
    
    According to equations~(\ref{eq: 14}), we can get~a linear equation with respect to~$m$ using the $n$-iteration estimate of $V_{\beta,m,n}(\omega_i)$ and $V_{\beta,m,n}(\mathcal{T}(\omega))$
    \begin{multline*}
        \left[1 + \beta\left(k_{n,0} - \sum_{i=1}^Kp_i(\omega)k_{n,i}\right)\right] m \\
        = \omega + \beta\left(\sum_{i=1}^Kp_i(\omega)a_{n,i} - a_{n,0}\right)
    \end{multline*}
    Given the belief state $\omega$, setting $\omega_\beta^*(m) = \omega$, if 
    \begin{equation*}
        1 + \beta\left(k_{n,0} - \sum_{i=1}^Kp_i(\omega)k_{n,i}\right) \neq 0,
    \end{equation*}
    we get the $n$-iteration approximated Whittle index as below
    \begin{equation*}
        \widehat{W}_n(\omega) = \left.\frac{\omega + \beta\left(\sum_{i=1}^Kp_i(\omega)a_{n,i} - a_{n,0}\right)}{1 + \beta\left(k_{n,0} - \sum_{i=1}^Kp_i(\omega)k_{n,i}\right)}\right|_{\omega^*_\beta(m) = \omega}.
    \end{equation*}
    Note that when $n = 0$, the $0$-iteration approximated Whittle index is the same as the imperfect Whittle index we computed before~(\ref{eq: 17}).
\end{IEEEproof}

\subsection{Algorithm and Complexity}
We summarize the above solution process into an algorithm called the Approximated Whittle Index (AWI) Policy. 
\begin{figure}[!t]
    \removelatexerror
    \begin{algorithm}[H]
        \caption{Approximated Whittle Index Policy}\label{alg:AWI}
        \begin{algorithmic}
            \REQUIRE $\beta \in(0,1), T \geq 1, N \geq 2, 1 \leq M<N, n_{iter} \geq 0$
            \REQUIRE $\omega_n(1),p_{11}^{(n)},p_{01}^{(n)},p_{i,1}^{(n)},p_{i,0}^{(n)},B_n\,(n=1,2,\cdots, N,\, i = 1,2,\cdots,K)$
            \FOR{$t = 1,2,\cdots,T$}
            \FOR{$n = 1,\cdots,N$}
            \STATE Set the threshold $\omega^*_\beta(m) = \omega_n(t)$ and try to compute the approximated Whittle index $\widehat{W}_{n_{iter}}(\omega_n(t))$
            \IF{$\widehat{W}_{n_{iter}}(\omega_n(t))$ exists}
            \STATE Set $W(\omega_n(t)) = \widehat{W}_{n_{iter}}(\omega_n(t))$
            \ELSE
            \STATE Set $W(\omega_n(t)) = \omega_n(t)B_n$
            \ENDIF
            \ENDFOR
            \STATE Choose the top $M$ channels with the largest approximated Whittle Indices $W(\omega_n(t))$
            \STATE Observe the selected $M$ channels and accrue reward $S_n(t)B_n$ from each active channel
            \FOR{$n=1,\cdots,N$}
            \STATE Update the belief state $\omega_n(t)$ according to (\ref{eq: 18})
            \ENDFOR
            \ENDFOR
        \end{algorithmic}
    \end{algorithm}
\end{figure}

\begin{theorem}\label{thm:complexity}
The complexity of Algorithm~\ref{alg:AWI} is $O(NTK^{n_{iter}})$, where~$N$ is the number of channels, $T$ the number of time steps, $K$ the number of CQI levels, and $n_{iter}$ the number of iteration steps in solving for the priority index. Furthermore, under Assumption~\ref{ass:threshold} and Assumption~\ref{ass:indexability} as $n_{iter}$ increases, the priority index converges to the exact Whittle index at a geometric rate.
\end{theorem}
\begin{IEEEproof}
The linear complexity in~$NT$ is obvious since arms are decoupled when computing their Whittle indices which are functions of their current belief states, respectively. From~\eqref{eq: 13}, the number of steps in calculating the value function at any belief state~$\omega$ is $O(1+K+K^2+\cdots+K^{n_{iter}})=O(K^{n_{iter}})$.  According to~\eqref{eq: 14}, the complexity of solving for the Whittle index of a single arm in each time slot is thus given by $O(K^{n_{iter}})$. If we look at~\eqref{eq: 14} again, the solution to subsidy~$m$ has the form $\frac{a}{b}$, where~$a$ is the expected total discounted reward under the active action in $V_{\beta,m}(\omega; a=1)$ minus that in $ V_{\beta,m}(\omega; a=0)$, and~$b$ is the expected total discounted time being passive in $ V_{\beta,m}(\omega; a=0)$ minus that in $V_{\beta,m}(\omega; a=1)$. When the iteration number is $n_{iter}$, the induced error in~$a$ or~$b$ is bounded by $\frac{2\beta^{n_{iter}}}{1-\beta}$. Then it is clear that the error in~$m$ also decreases geometrically with $n_{iter}$.
\end{IEEEproof}

If Assumption~\ref{ass:threshold} or Assumption~\ref{ass:indexability}, or both of them are not verified, the optimal single-armed policy is not guaranteed to have a threshold structure and the original single-armed problem is not guaranteed to be indexable. Consequently, the priority index calculated by our algorithm may not converge to the true Whittle index. However, to implement Algorithm~\ref{alg:AWI}, we do not need to satisfy the above assumptions. According to the continuity of the value function, for any given~$\omega$, there exists at least one subsidy~$m$ such that~$\omega$ is a threshold under~$m$. However, the exact system of linear equations may not have a solution when the threshold is not unique, i.e. even relaxed indexability does not hold~\cite{liu_rmab_2025}. In this case, if the system of linear equations has no solution, we fall back to the myopic index. In our numerical experiments, 20,768 instances were able to solve for the index and achieve good performance.

\section{Simulation Results}
\label{sec: simulation results}
In this section, we evaluate the proposed approximate Whittle index policy through a set of numerical experiments. The experiments are designed to examine its effectiveness, robustness, parameter sensitivity, scalability, and computational cost. The notation used throughout the numerical studies is summarized in Table~\ref{tab:notations}, while the threshold optimality and indexability conditions of the test instances are classified in Table~\ref{tab:threshold_indexability_cases}. The main experimental results are reported in Figs.~\ref{fig:baseline_comparison}--\ref{fig:runtime}. In practice, we terminate the iteration when the ranking of the Whittle indices remains unchanged between two consecutive iterations. Across all 20,768 randomly generated instances considered in our experiments, the index ranking stabilizes within at most three iterations, indicating that AWI-3 is sufficient. AWI-3 is therefore omitted from the corresponding figures when it produces the same index ranking, selects the same channels, and achieves the same performance as AWI-2.

For each channel $n$, the transition probabilities
$\bigl(p_{01}^{(n)},p_{11}^{(n)}\bigr)$ are randomly generated. For each stationary Gilbert--Elliott channel, we assume that the channel has evolved for a sufficiently long time before the beginning of the simulation, which is defined as time slot $t=0$. Therefore, the initial state is sampled from the stationary distribution.

At each time slot, a policy computes a score for every channel and selects the active channels. The immediate reward is then collected, and a CQI observation is randomly generated for each active channel according to
$\{p_{i,0}^n:i=1,\ldots,K\}$ and
$\{p_{i,1}^n:i=1,\ldots,K\}$, according to the current states of channels. The channel states evolve according to their Gilbert-Elliott channel model. Finally, we update the belief states in accordance to (\ref{eq: 18}) for the next round of decision-making.

For a finite horizon of $T$ time slots, performance is measured by the average discounted return over independent Monte Carlo runs, where
\begin{equation}
    G_\pi(T)=\sum_{t=0}^{T-1}\beta^t R_\pi(t),
\end{equation}
and $R_\pi(t)$ denotes the reward obtained by policy $\pi$ at time slot $t$.

For detailed experimental setting , please see Table.~\ref{tab:experiment_statistics}

\subsection{Comparison with Baselines and Exact Dynamic Programming}
\label{subsec:baseline_exact}

\paragraph{Baseline implementations}
The first baseline is the Myopic policy, which selects the
$M$ channels with the largest immediate expected rewards:
\begin{equation}
    \pi_{\mathrm{myopic}}(\boldsymbol{\omega})
    \in
    \arg\max_{\mathcal{A}\subseteq\{1,\ldots,N\}:\,|\mathcal{A}|=M}
    \sum_{n\in\mathcal{A}}\omega_n B_n.
\end{equation}

The second baseline, denoted by AWI0, uses the approximate
Whittle-index computation proposed in
\cite{elmaghraby_femtocell_2018}.

The third baseline is Myopic-Rollout, implemented as a three-stage finite-horizon lookahead based on the dynamic programming principle~\cite{bertsekas2020rollout}. Starting from the current belief state, all feasible actions and their reachable observation outcomes are evaluated recursively over three decision stages. At each intermediate stage, the action maximizing the expected discounted return over the remaining horizon is selected. At the last stage, since no future reward is considered, the optimal decision reduces to the Myopic action that maximizes the immediate expected reward. The current action with the largest three-stage expected discounted return is then selected.

The fourth baseline is value function approximation (VFA),
implemented as a finite-horizon, model-based adaptation of the
tree-based batch reinforcement-learning framework in
\cite{ernst2005tree}. The continuation value is approximated in
the joint belief space using fitted Bellman targets. During
online scheduling, VFA performs a two-step Bellman lookahead
and then uses the fitted continuation value to evaluate each
feasible current action. Unlike the Myopic policy, VFA therefore
accounts for both future belief evolution and subsequent action
reoptimization.

\begin{figure*}[!t]
\centering
\subfloat[$\beta=0.5$.]{
    \includegraphics[width=0.4\textwidth]{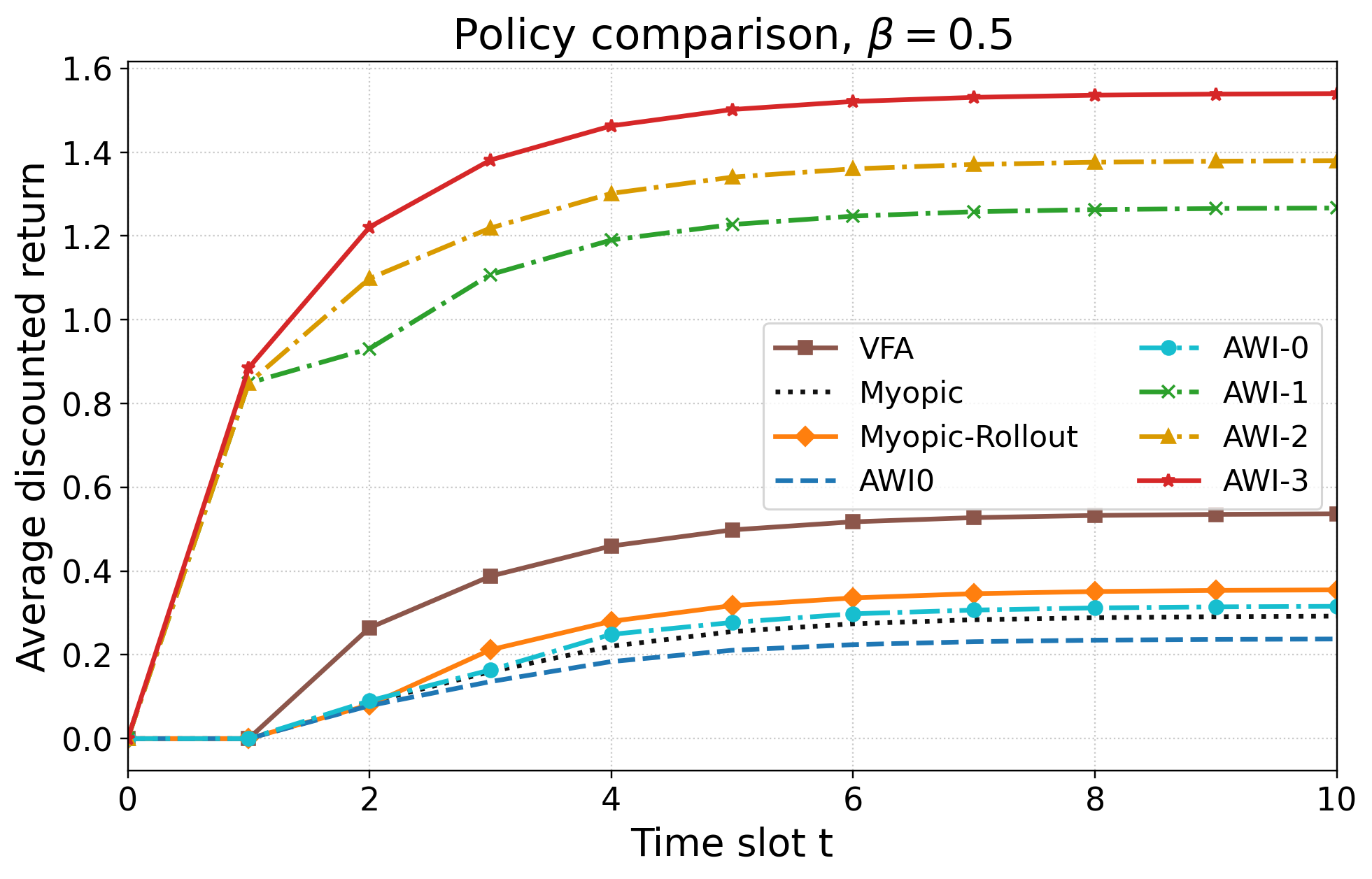}
    \label{fig:baseline_beta05}
}
\hfil
\subfloat[$\beta=0.9$.]{
    \includegraphics[width=0.4\textwidth]{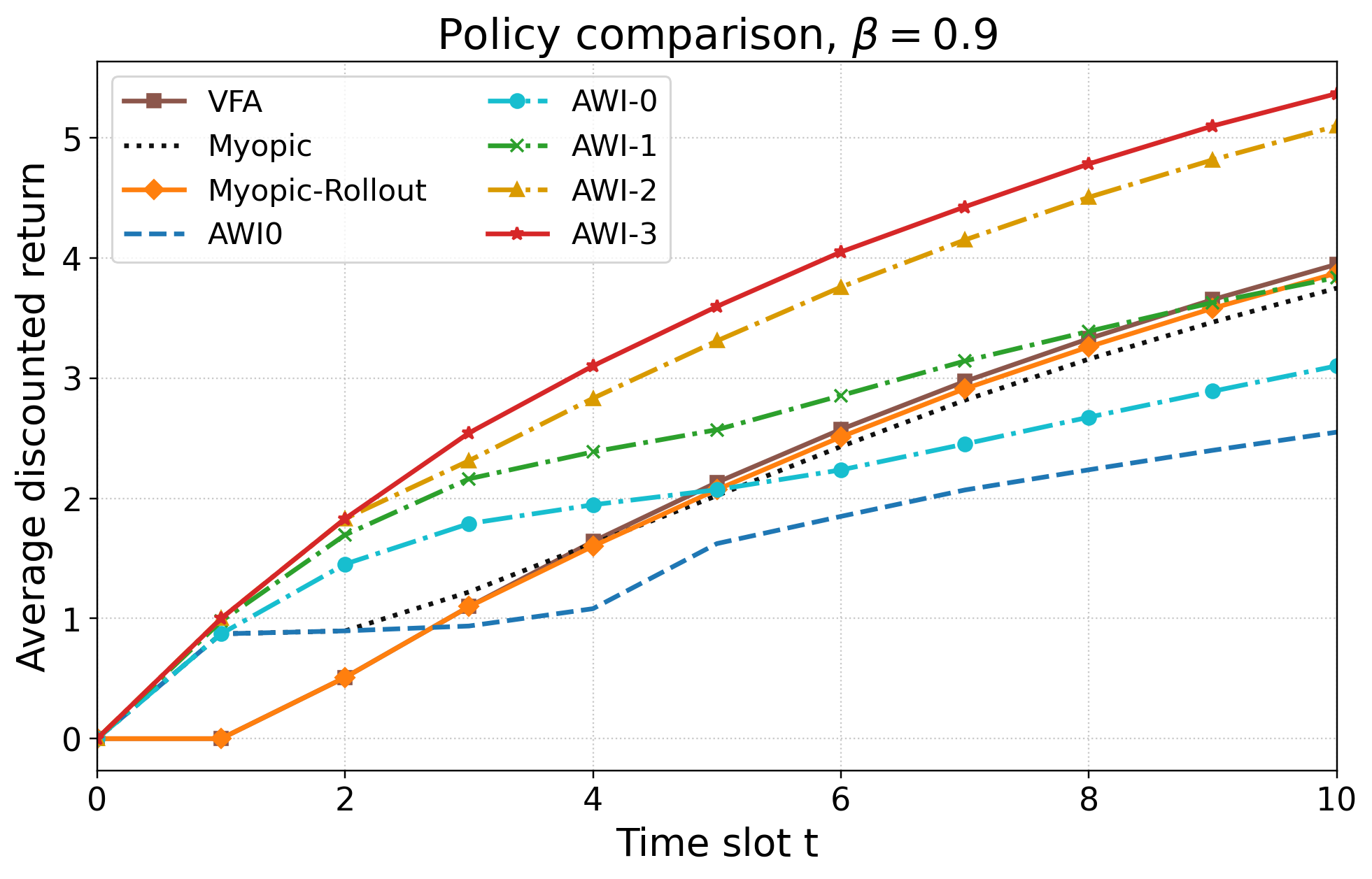}
    \label{fig:baseline_beta09}
}
\caption{Average discounted return of the proposed AWI policies and the baseline policies under different discount factors.}
\label{fig:baseline_comparison}
\end{figure*}

Fig.~\ref{fig:baseline_comparison} shows that increasing the AWI refinement iteration consistently improves the average discounted return in both instances. AWI-0 already improves upon the existing AWI0 approximation, while AWI-1, AWI-2, and AWI-3 progressively incorporate more future belief. Among all compared policies, AWI-3 achieves the highest return under both $\beta=0.5$ and $\beta=0.9$. In particular, it outperforms the Myopic-Rollout and the VFA baseline, while requiring substantially less online computation, as further demonstrated by the runtime experiments according to Fig.~\ref{fig:runtime_arms}

To directly evaluate the AWI policy, we further construct a finite-horizon instance for which the joint belief-state dynamic program can be solved exactly. This benchmark has $N=8$, $M=1$, $K=2$, $T=6$, and $\beta=0.9$. Exact DP exhaustively enumerates all feasible channel selections and all possible CQI observations over the complete reachable joint belief tree. The optimal finite-horizon return is computed by backward induction, with a zero terminal value after the final decision epoch.

\begin{figure*}[!t]
\centering
{
    \includegraphics[width=0.45\textwidth]{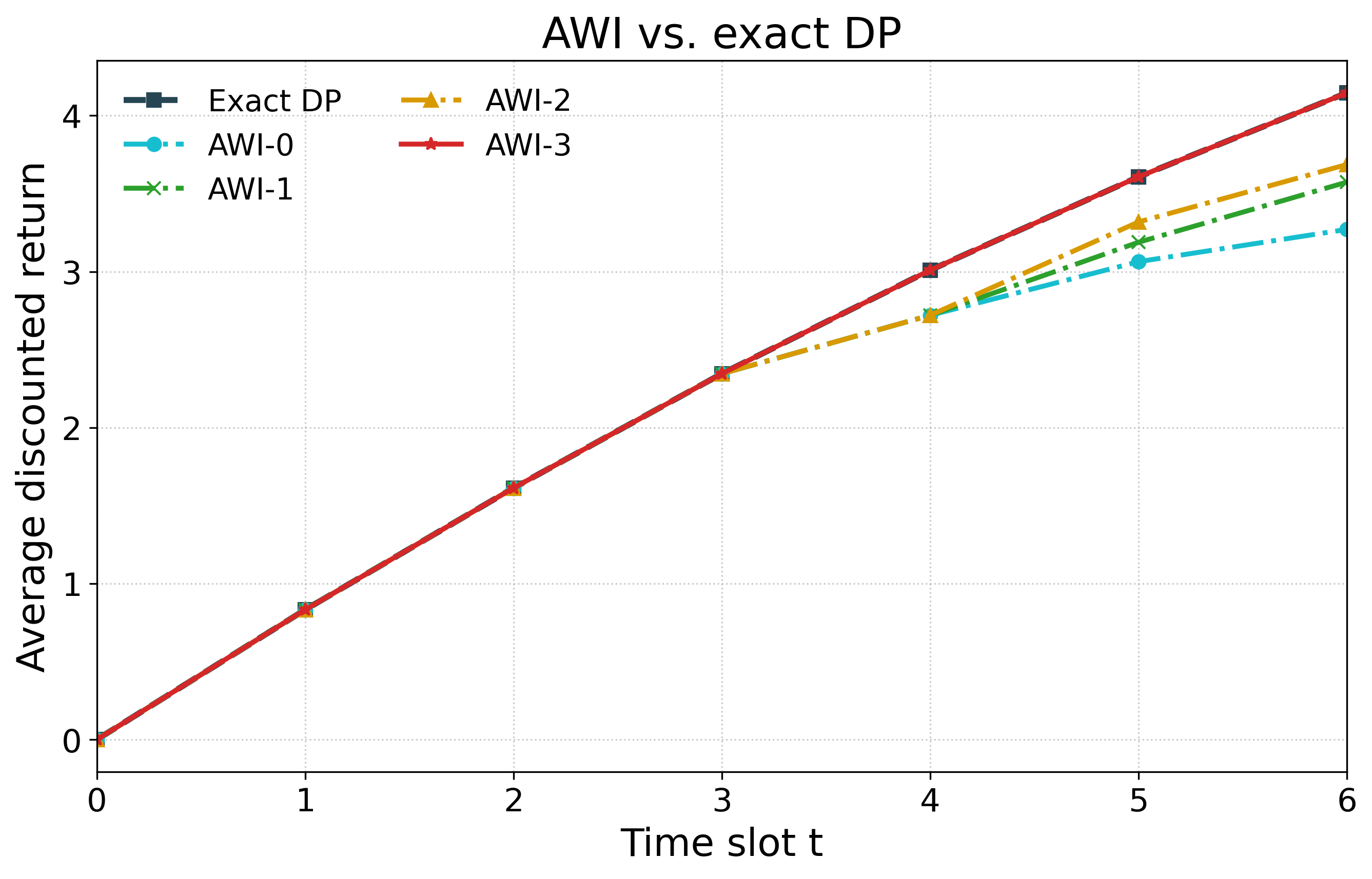}
    \label{fig:exact_time}
}
\hfil
{
    \includegraphics[width=0.45\textwidth]{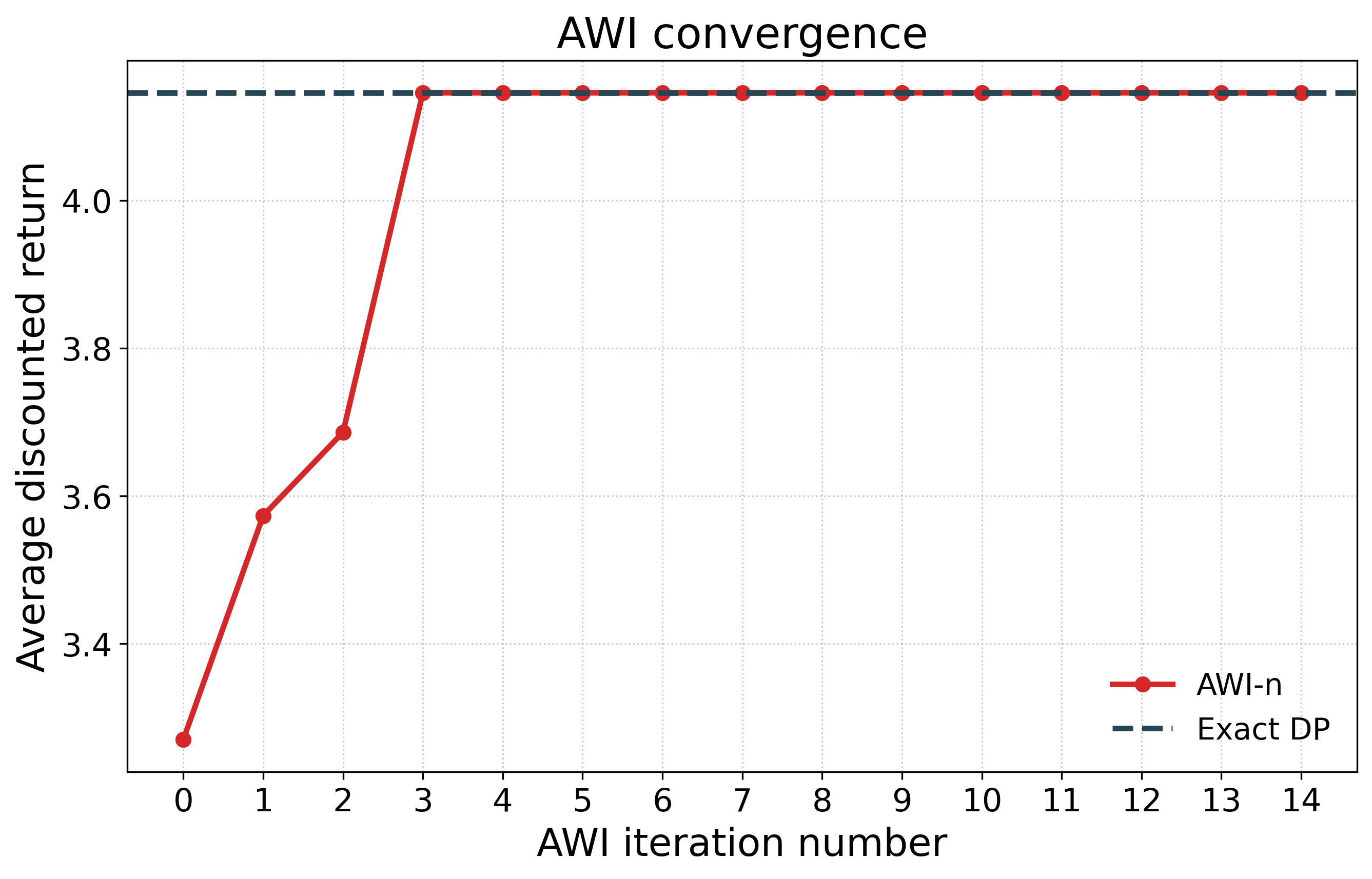}
    \label{fig:exact_depth}
}
\caption{Comparison between the AWI policies and Exact DP on the finite-horizon exact benchmark.}
\label{fig:exact_comparison}
\end{figure*}

As shown in Fig.~\ref{fig:exact_comparison}, the exact expected discounted returns of AWI-0, AWI-1, and AWI-2 are $3.2704$, $3.5734$, and $3.6862$, respectively. AWI-3 attains $4.1455$, which is identical to the Exact DP value. The complete return trajectory of AWI-3 also the same as Exact DP.

Fig.~\ref{fig:exact_comparison} further shows that the return increases with the recursion depth and reaches the Exact DP value at $n=3$.

Rollout-myopic and VFA policies are far more computationally expensive than the AWI policy without providing a clear performance advantage. We therefore omit them from the remaining experiments to maintain figure clarity. The myopic policy alone serves as a consolidated proxy for these simple baselines. Exact dynamic programming becomes computationally intractable for large-scale instances and is therefore used only in small scale experiments as an upper-bound benchmark.

\subsection{Large-System Performance}

We next evaluate the proposed method on two large-system instances. The compared policies are Myopic,
AWI0, and AWI-$n$ for $n=0,\ldots,3$. 

\begin{figure*}[!t]
\centering
\subfloat[$\beta=0.5$.]{
    \includegraphics[width=0.4\textwidth]
    {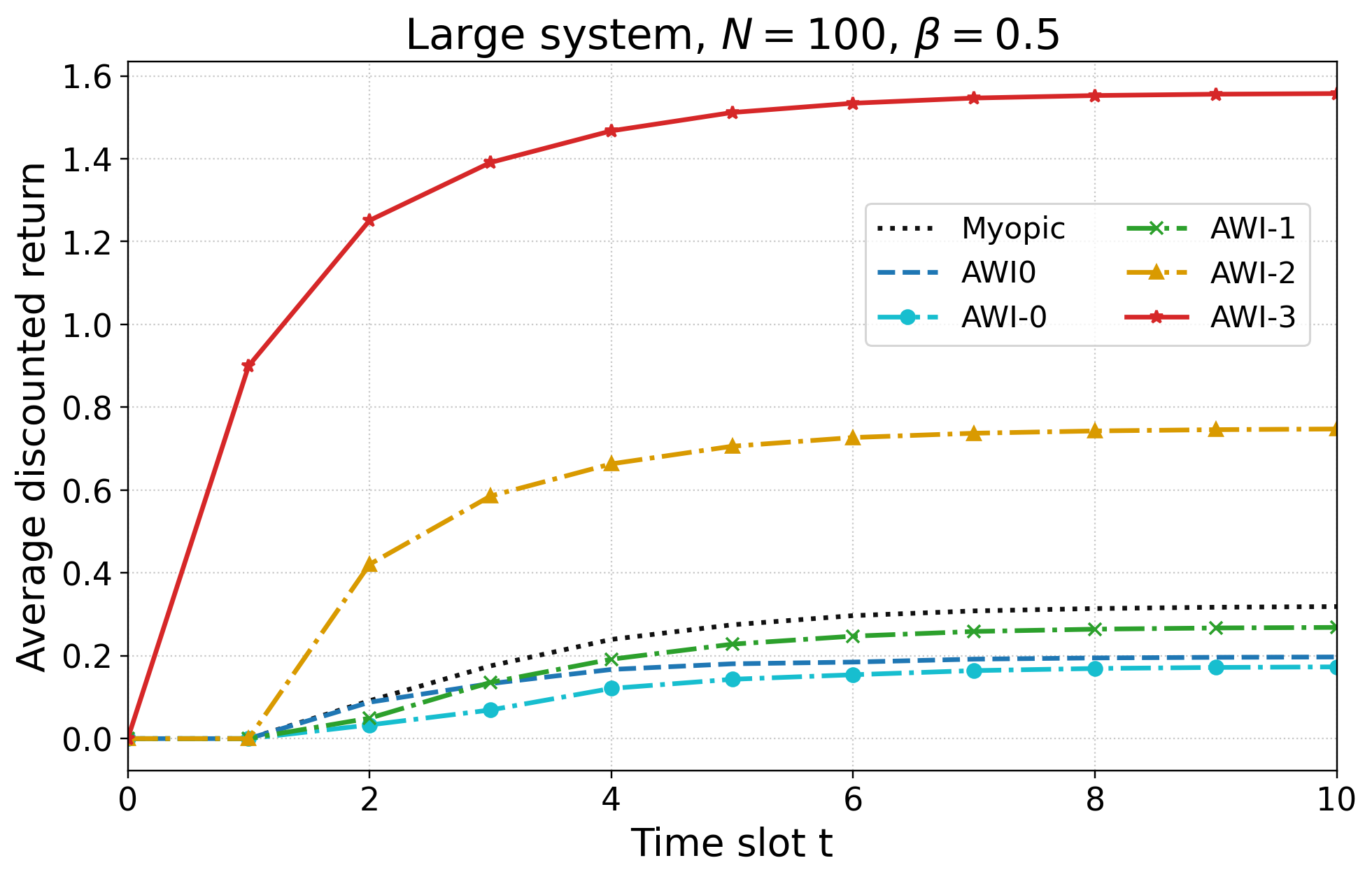}
    \label{fig:large_beta05}
}
\hfil
\subfloat[$\beta=0.9$.]{
    \includegraphics[width=0.4\textwidth]
    {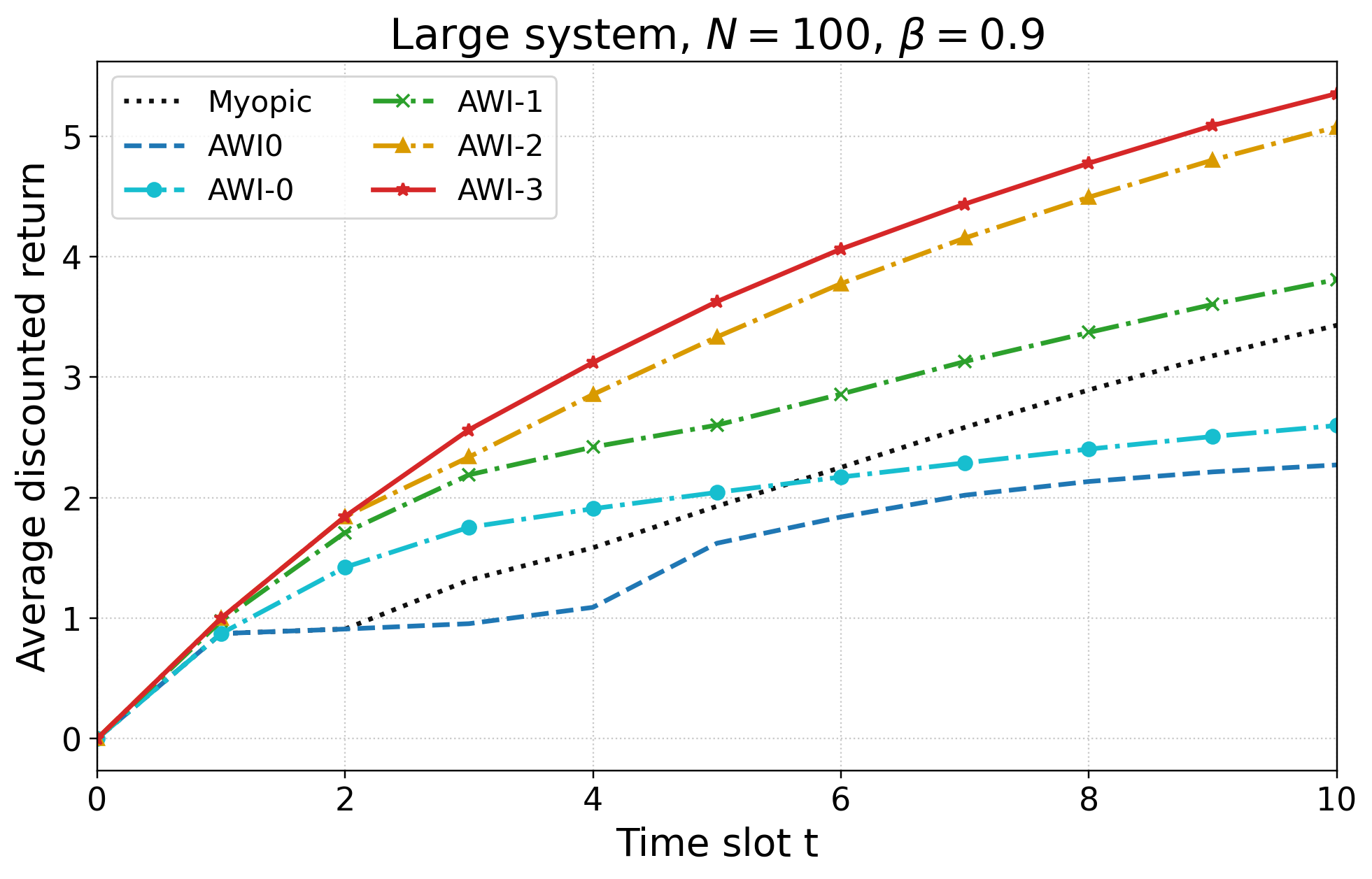}
    \label{fig:large_beta09}
}
\caption{large-system experiments}
\label{fig:large_system}
\end{figure*}

As shown in Fig.~\ref{fig:large_system}, AWI-3 achieves the highest
average discounted return in both large system instances. These results demonstrate that the performance advantage of the AWI policy is retained in the tested systems with $100$ channels.

\subsection{Non-Stationary Channel Dynamics}
We further consider two non-stationary channel environments  where the transition probabilities change during the simulation horizon. The change point is marked in the figures.

\begin{figure*}[!t]
\centering
\subfloat[$\beta=0.5$.]{
    \includegraphics[width=0.4\textwidth]{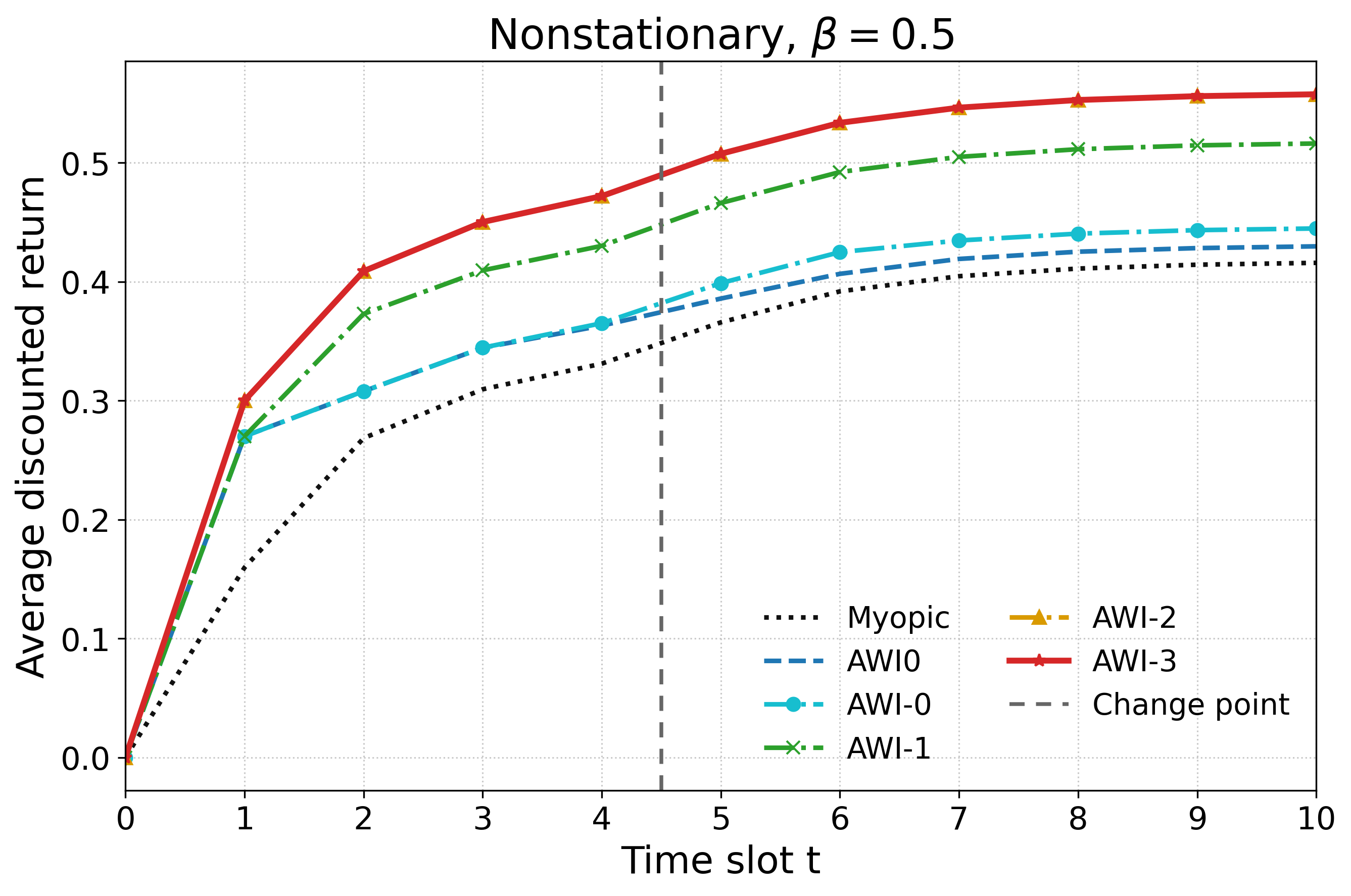}
    \label{fig:nonstationary_beta05}
}
\hfil
\subfloat[$\beta=0.9$.]{
    \includegraphics[width=0.4\textwidth]{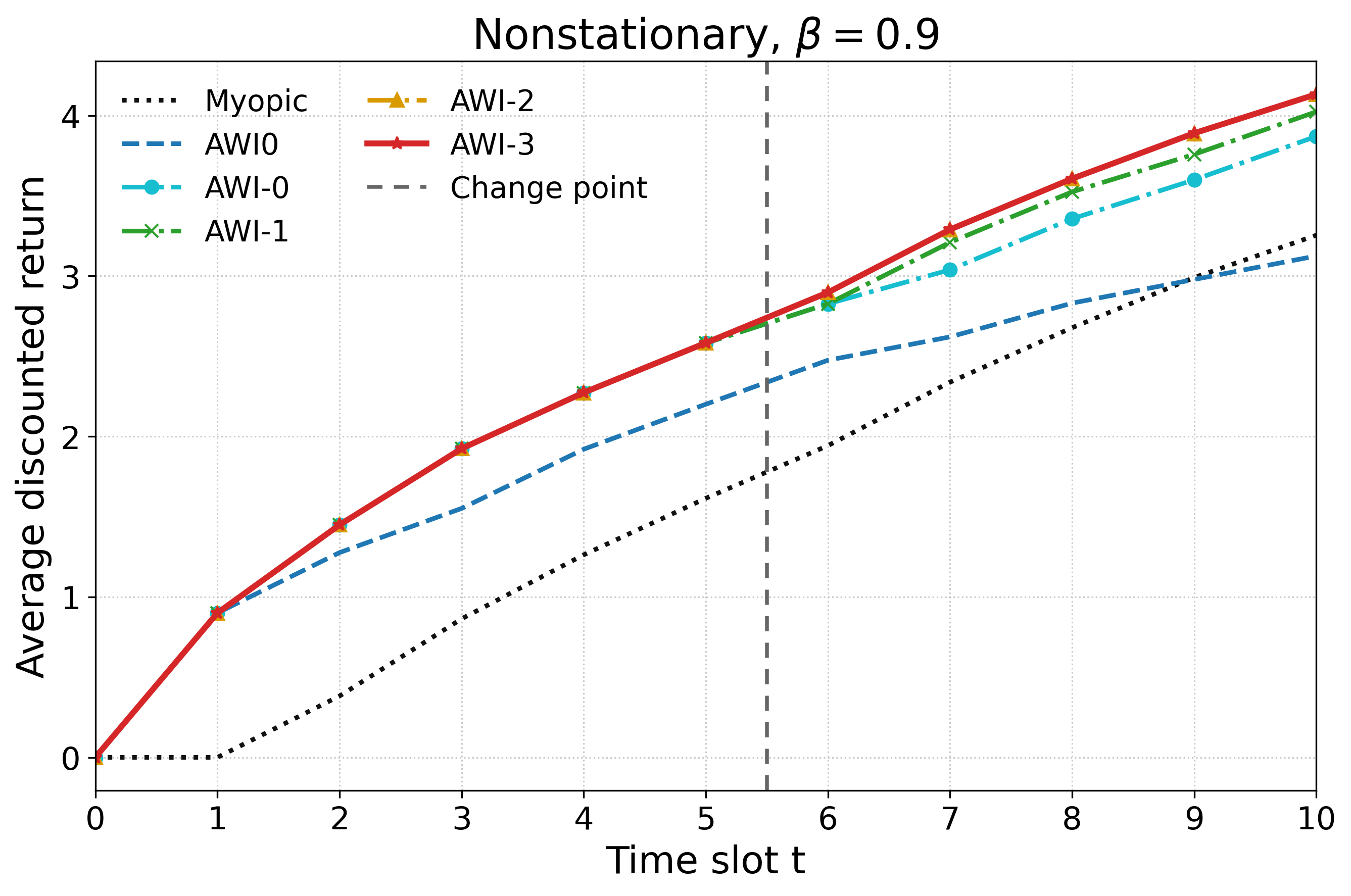}
    \label{fig:nonstationary_beta09}
}
\caption{Performance under non-stationary channel dynamics.}
\label{fig:nonstationary}
\end{figure*}

As shown in Fig.~\ref{fig:nonstationary}, the proposed AWI policies continue to perform well after the parameters changed. This indicates that the proposed method is not restricted to a fixed stationary channel setting and remains robust when the channel dynamics vary over time.

\subsection{Sensitivity to CQI}

We further evaluate the sensitivity of the proposed method to the number of CQI levels. Different CQI configurations lead to different observation structures and hence different belief updates.

\begin{figure*}[!t]
\centering
\subfloat[CQI setting 1.]{
    \includegraphics[width=0.4\textwidth]{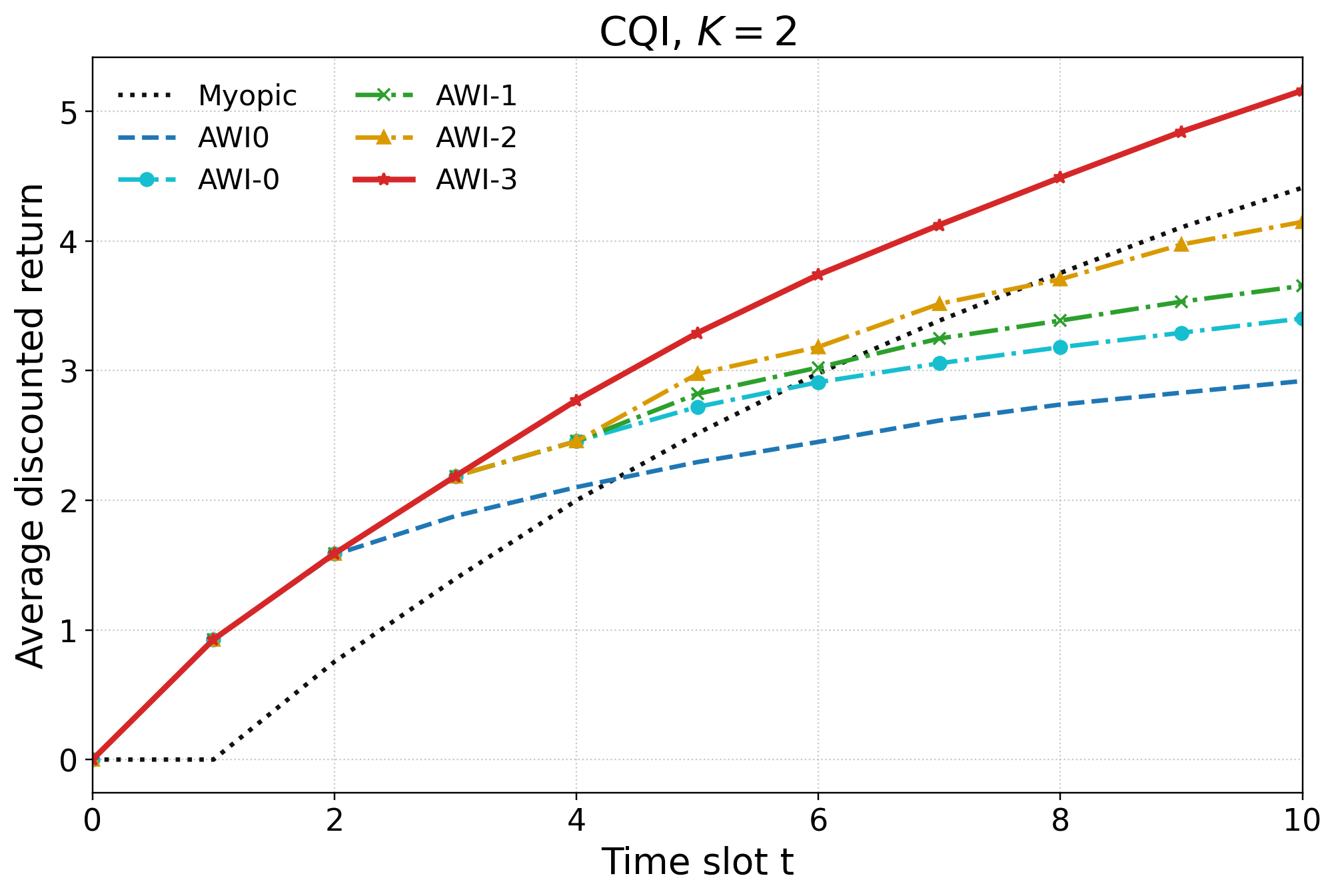}
    \label{fig:cqi_setting1}
}
\hfil
\subfloat[CQI setting 2.]{
    \includegraphics[width=0.4\textwidth]{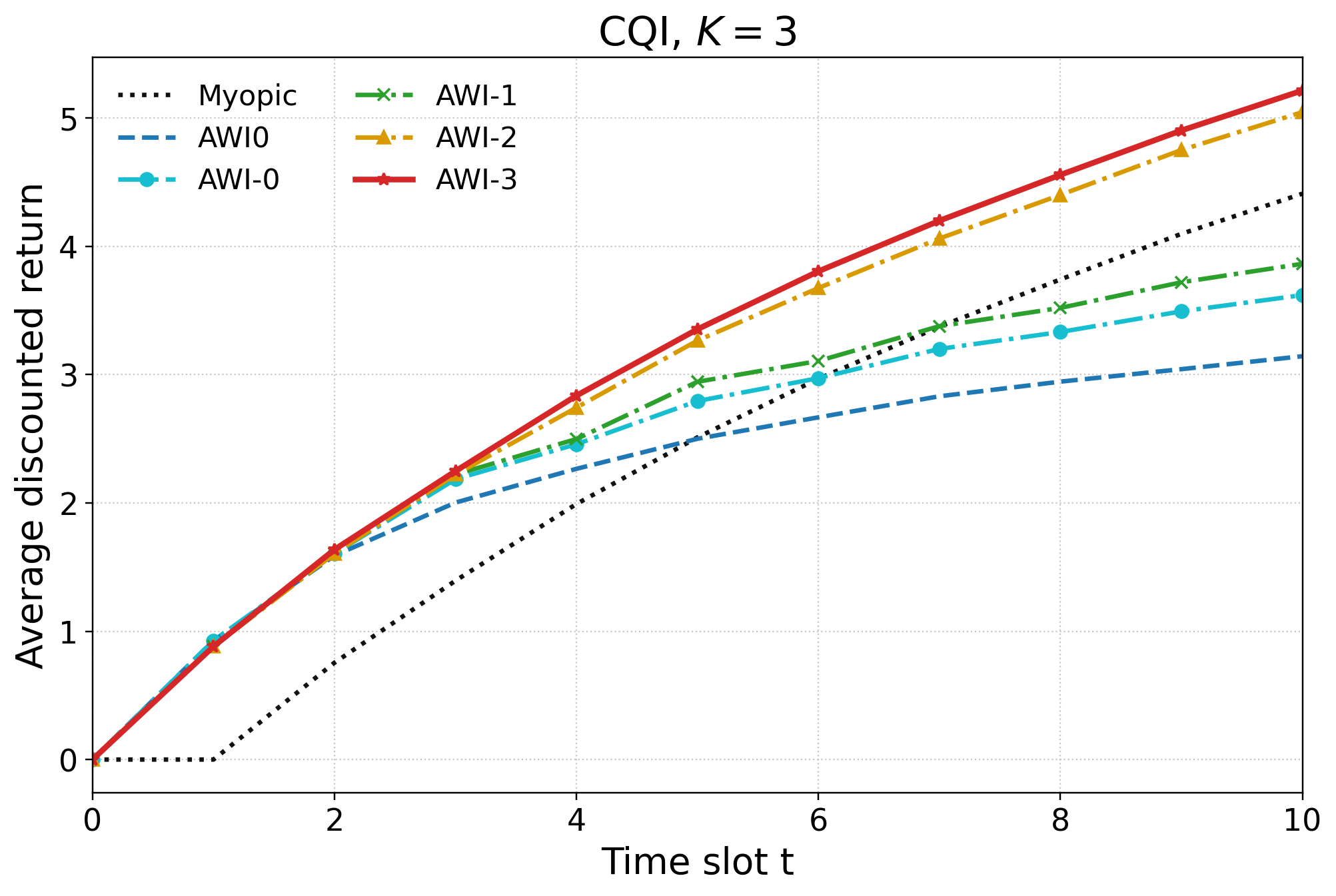}
    \label{fig:cqi_setting2}
}

\vspace{0.4em}

\subfloat[CQI setting 3.]{
    \includegraphics[width=0.4\textwidth]{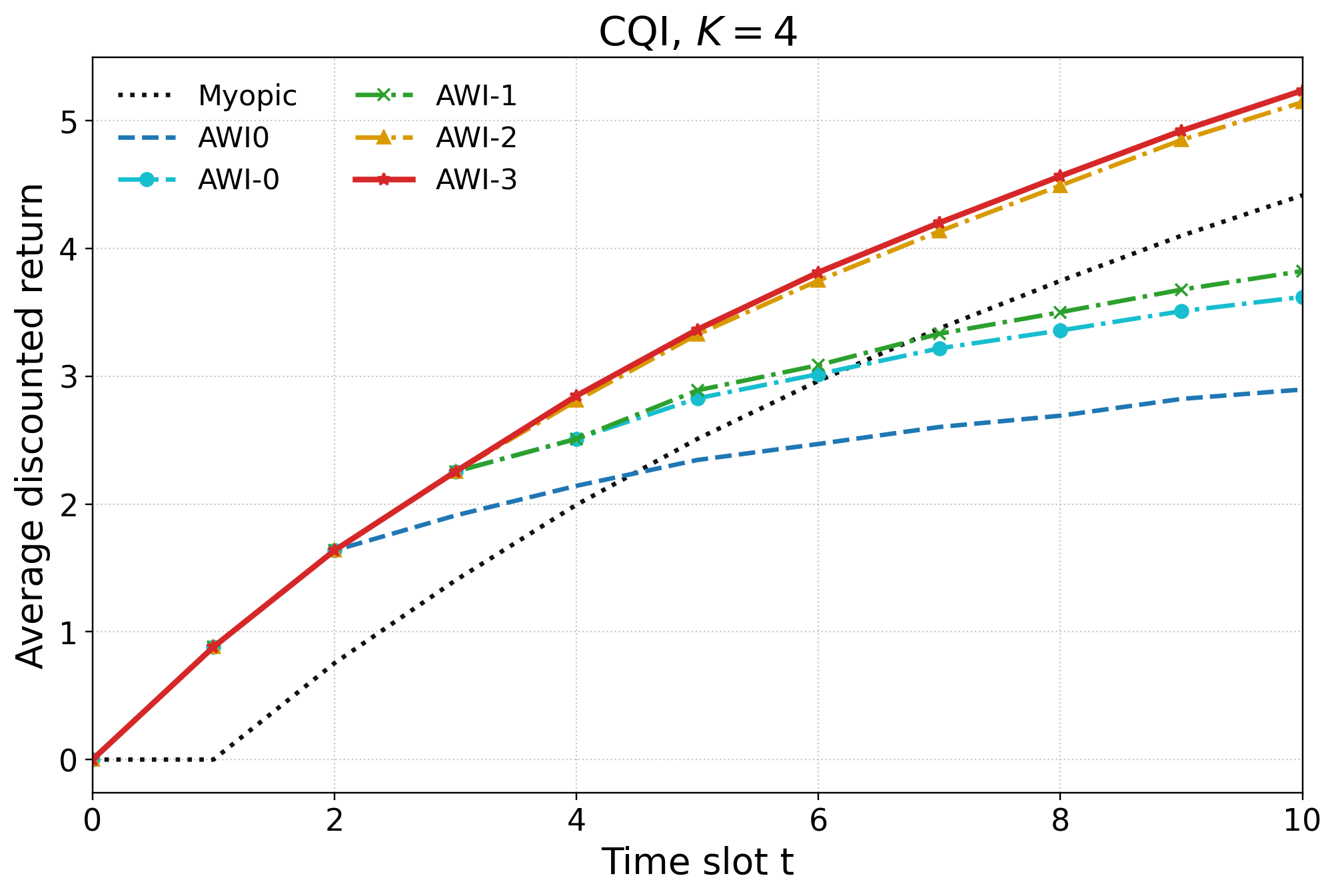}
    \label{fig:cqi_setting3}
}
\hfil
\subfloat[CQI setting 4.]{
    \includegraphics[width=0.4\textwidth]{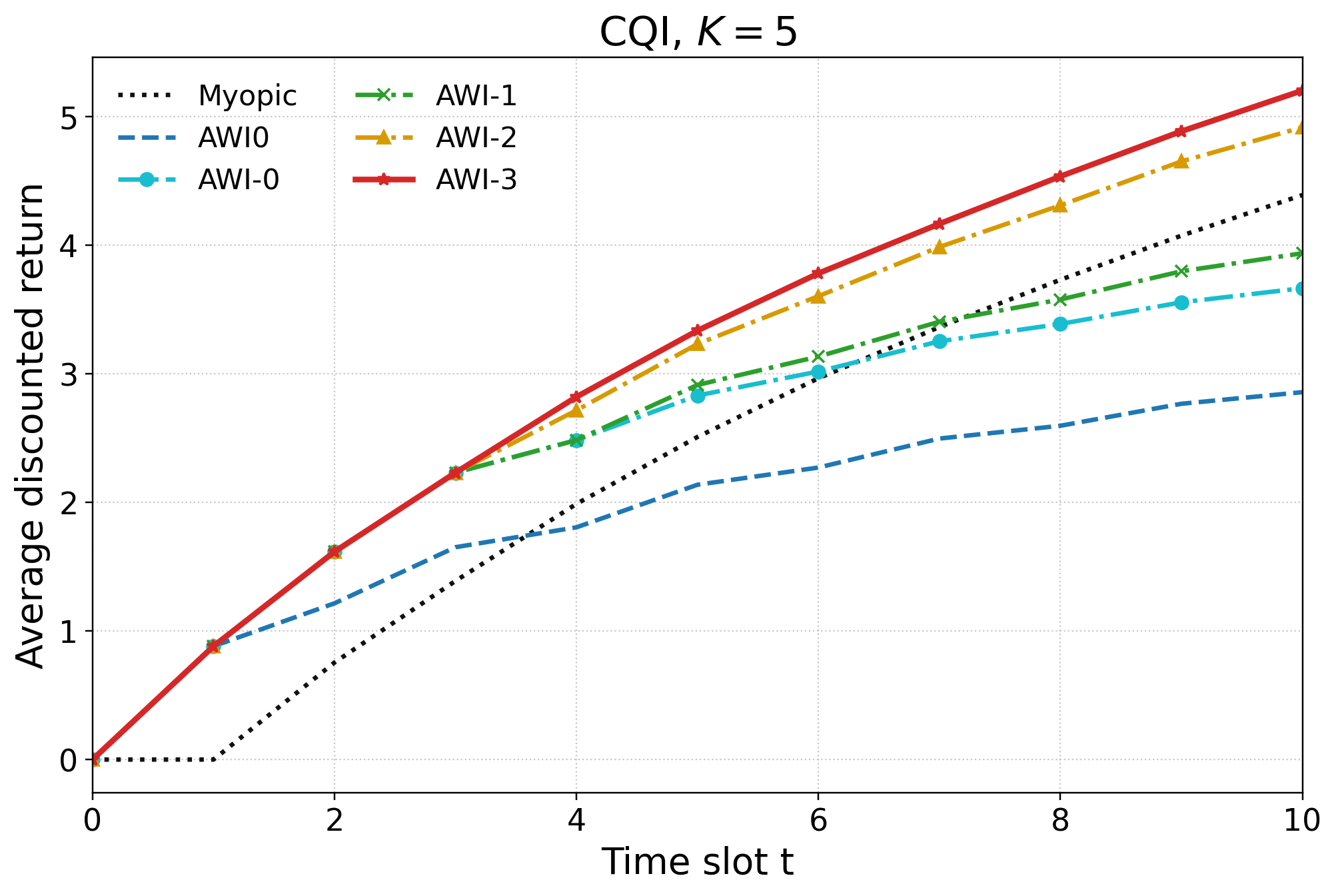}
    \label{fig:cqi_setting4}
}
\caption{Sensitivity analysis with respect to CQI configurations.}
\label{fig:cqi_sensitivity}
\end{figure*}

Fig.~\ref{fig:cqi_sensitivity} shows that the proposed AWI policies consistently outperform the Myopic policy under different CQI configurations. This suggests that the performance gain is not caused by one particular observation setting, but is preserved when the CQI model changes.

\subsection{Simulation with $\beta=1$}

Our numerical experiments showed that the Whittle index converges as the discount factor $\beta$ approaches one. 

\begin{figure*}[!t]
\centering
\includegraphics[width=0.4\textwidth]{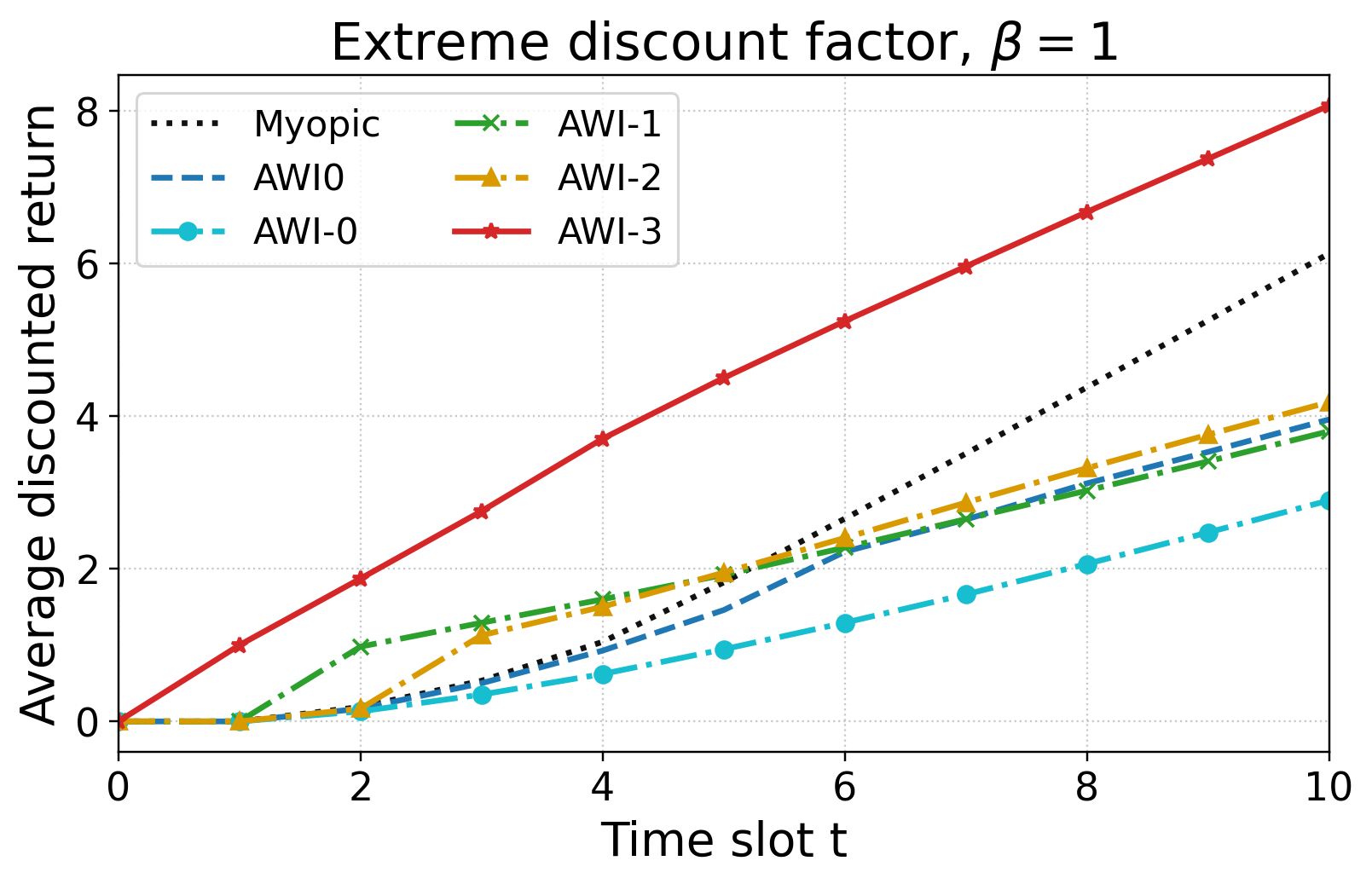}
\caption{Performance under an extreme high-discount factor $\beta=1$.}
\label{fig:beta09999}
\end{figure*}

Fig.~\ref{fig:beta09999} shows that the proposed AWI policies remain effective. In this case, the Myopic policy becomes less competitive because it ignores the long-term impact of current actions on future belief states. 

\subsection{Runtime Analysis}

Finally, we evaluate the computational cost of the proposed iteration. According to Theorem~\ref{thm:complexity}, the total
complexity of Algorithm~\ref{alg:AWI} over $T$ time slots is
$O(NTK^{n_{\mathrm{iter}}})$, corresponding to
$O(NK^{n_{\mathrm{iter}}})$ per scheduling decision. Therefore, for
fixed $K$ and $n_{\mathrm{iter}}$, the computational cost grows linearly with the number of channels, whereas increasing the refinement depth introduces the principal additional cost.

\begin{figure*}[!t]
\centering
\subfloat[Runtime versus the number of channels.]{
    \includegraphics[width=0.4\textwidth]
    {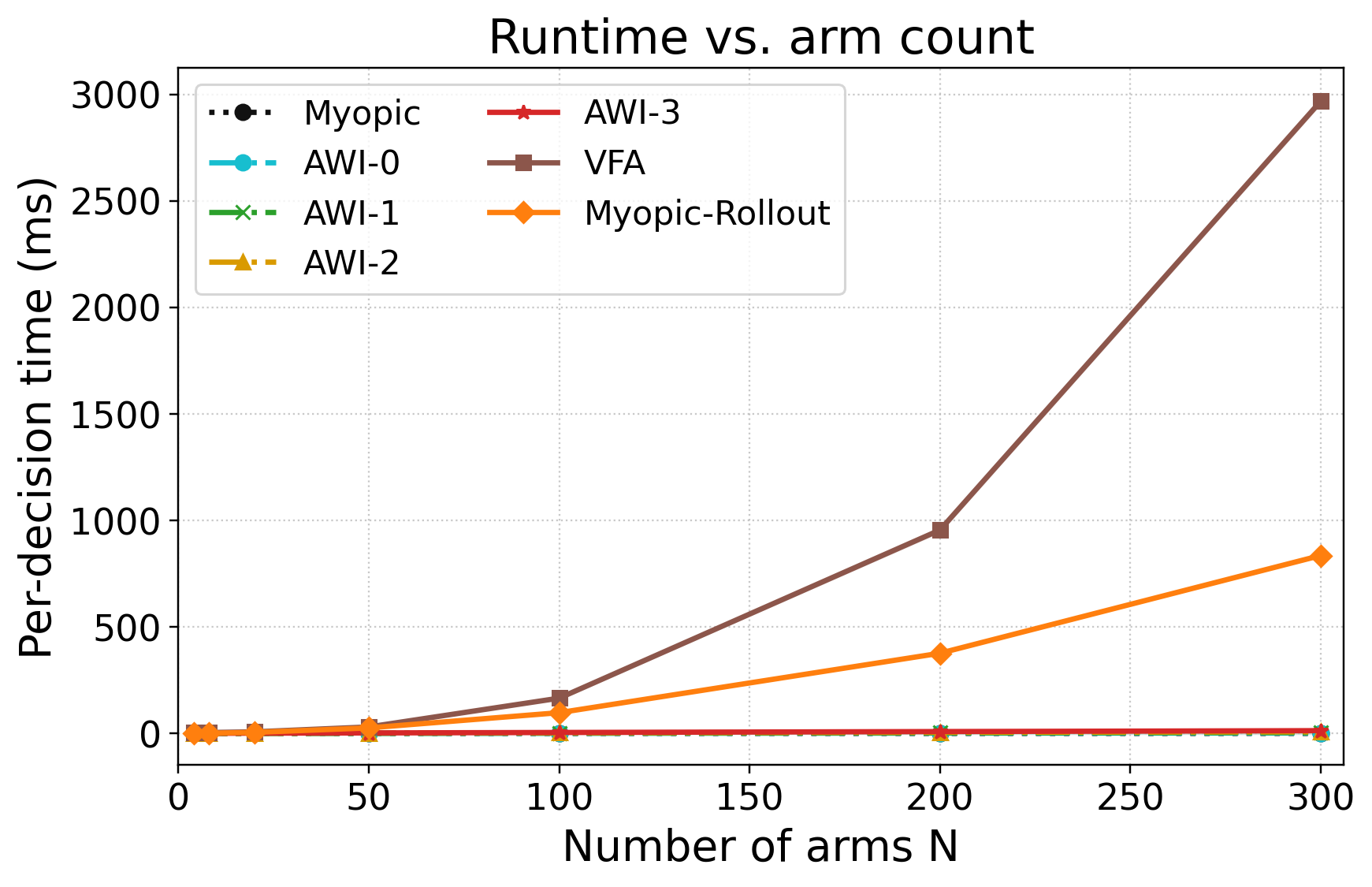}
    \label{fig:runtime_arms}
}
\hfil
\subfloat[AWI index computation time versus depth at $N=200$.]{
    \includegraphics[width=0.4\textwidth]
    {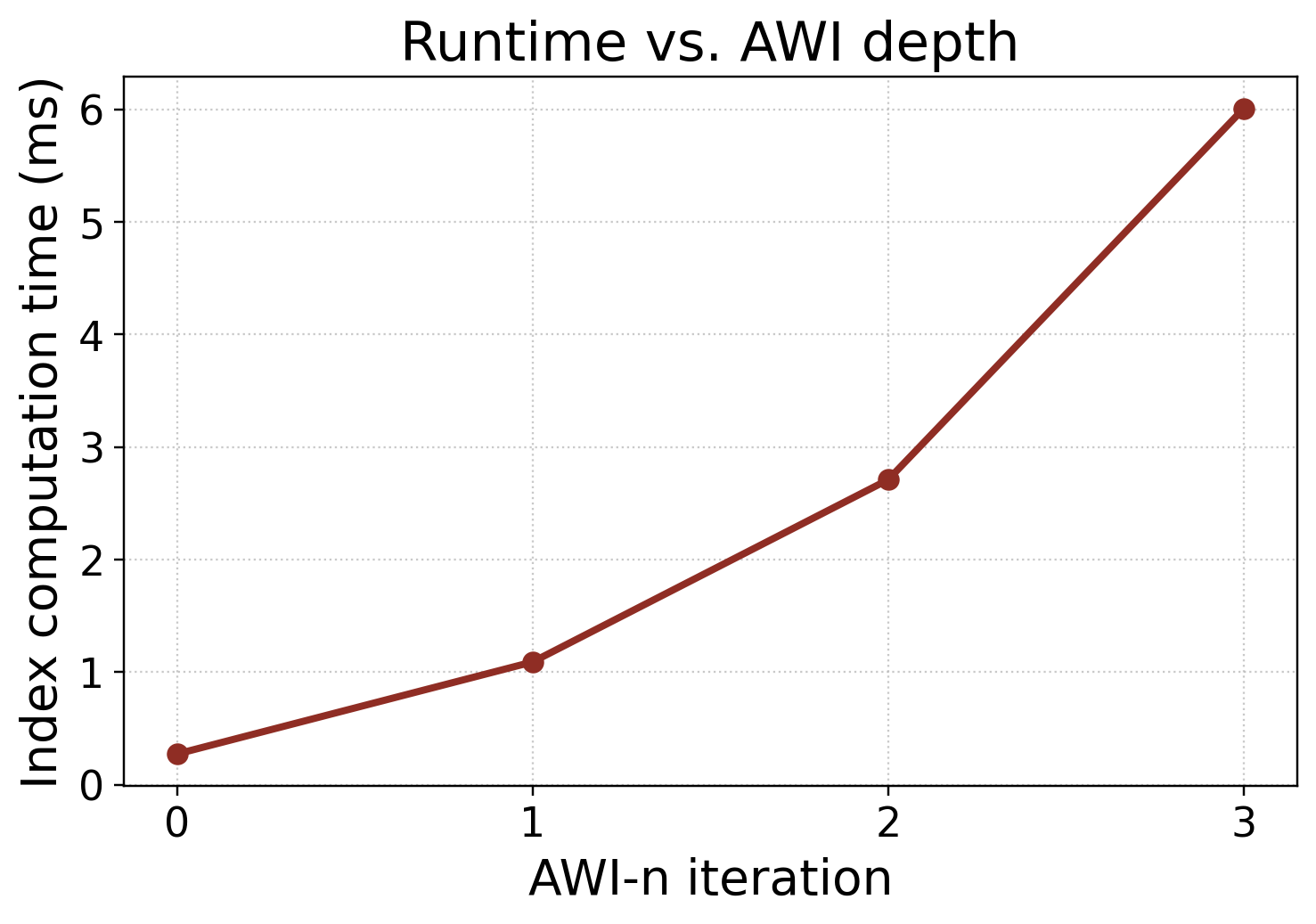}
    \label{fig:runtime_iterations}
}
\caption{Online runtime and scalability comparison.}
\label{fig:runtime}
\end{figure*}

For the experiment setting, the number of channels is
varied over $N\in\{4,8,20,50,100,200,300\}$. As shown in Fig. \ref{fig:runtime_arms}, the online cost of each AWI and Myopic policy grows approximately linearly with $N$ and remains substantially lower than those of Myopic-Rollout and VFA. In contrast, VFA and Myopic-Rollout policies exhibit much steeper, superlinear growth, reaching approximately $2967.42$~ms and $835.02$~ms per decision, respectively, at N=300. These results demonstrate the favorable computational scalability of the our algorithm.  Fig.~\ref{fig:runtime_iterations} further shows that the AWI computation increases with the refinement depth, illustrating the tradeoff between index refinement and computational cost.

All runtime experiments were conducted on a MacBook Pro equipped with
an Apple M4 Pro processor (12 CPU cores: 8 performance and 4 efficiency
cores) and 24~GB of memory, running macOS 26.5.2 and Python 3.14.0.

Table~\ref{tab:experiment_statistics} summarizes the average discounted returns and associated statistical measures of all compared policies across the numerical experiments. Table~\ref{tab:threshold_indexability_cases}  classifies the experimental settings according to whether the conditions for threshold optimality and indexability are satisfied. Table~\ref{tab:comparison}  compares the proposed method with the binary-feedback RMAB approach~\cite{liu_low-complexity_2022} and our preliminary conference work~\cite{elmaghraby_femtocell_2018} in terms of model structure, theoretical analysis, algorithm design, and numerical evaluation. Table~\ref{tab:baseline_parameters}  lists the transition probability, channel bandwidths, and CQI probability used in the two baseline systems.

\section{Conclusion}
\label{sec: conclusion}
In this work, we employed the RMAB framework to model the DSA problem where CQI can be observed through cognitive capabilities of the sub-network. We established the indexability of the relaxed single-armed bandit problem and demonstrated the optimality of the threshold policy under certain conditions. Furthermore, we analyzed the value functions with tight approximations to obtain the Whittle Index Policy, a heuristic and low-complexity solution for such complex sequential scheduling problems. Finally, we substantiated the effectiveness and robustness of the proposed algorithm through extensive numerical simulations including the scenario where the conditions for indexability do not hold. 

In future research, it is worthwhile to investigate the impact of the iteration step size on the performance of the proposed policy and optimize it according to different problem settings. Additionally, the integration of deep learning and neural networks could be explored to enhance the computational efficiency in solving for the Whittle index. For example, we could train a neural network to search for a fast convergence path to the exact Whittle index by leveraging the rich structures of the value functions. 

Assuming unknown transition and observation probabilities substantially generalizes the problem by introducing a broader research scope. Specifically, estimating unknown system parameters falls into the non-Bayesian (frequentist) framework\cite{liu2026extended} of bandit theory, where the system performance is usually measured by regret (the performance loss due to unknown parameters compared with the optimal performance when the parameters are known). In general, these probabilities can be estimated by counting the frequencies of occurrences of the corresponding events in the observation history. However, ensuring sufficient data for rapid convergence requires an additional layer of exploration and fundamentally complicates the decision-making process of channel selection. Without sufficient exploration, the estimation errors might be too large to obtain an accurate Whittle index computation due to the heavy dependence of the belief transitions and value functions on those probabilities. This constitutes a challenging but promising direction for future research.

Future work also includes multi-secondary-user interference by incorporating game theory and decentralized bandit; nonstationary Markov processes by investigating into quickest detection theory and long-memory time series analysis; multi-state channel model by high-dimensional functional analysis; correlation among channels by joint design of index functions; unknown transition probabilities by incorporating results from the non-Bayesian class of bandit problems; etc.

\noindent\textbf{Acknowledgment} We thank the anonymous reviewers for the valuable comments that helped improve this paper.

\bibliographystyle{IEEEtran}
\bibliography{IEEEabrv,refs}

\bigskip

\noindent\textbf{Keqin Liu} received his PhD degree from the University of California at Davis in 2010 and joined industry as a software engineer in 2012 after 2 years of postdoc experience. After 8 years of working in the Silicon Valley, Dr. Liu still has a lot passion in academic research and decided to return to academia in 2020 when working in ASML and subsequently taught in Nanjing University until joining Xi'an Jiaotong-Liverpool University in 2024. Dr. Liu's research interests include the development of modern math theory for AI technologies and the application of pure theory to solving practical problems in machine learning, e.g., multi-armed bandit in reinforcement learning.

\bigskip

\noindent\textbf{Qizhen Jia} received the B.Sc. degree in financial mathematics from Xi'an Jiaotong-Liverpool University, China, in 2025. He is currently pursuing the Ph.D. degree with the School of Mathematics and Physics, Xi'an Jiaotong-Liverpool University, China. His research interests include restless multi-armed bandits, partially observable Markov decision processes, index policies, approximate dynamic programming, reinforcement learning, and stochastic optimization, with applications to wireless communications and dynamic spectrum access.

\bigskip

\noindent\textbf{Yiying Zhang} received the B.S. degree in mathematics from Nanjing University, Nanjing, China, in 2022 and the M.S. degree in mathematics from Nanjing University in 2025. During her M.S. studies, her research interests include reinforcement learning and multi-armed bandits. Currently, she is a Recommendation Algorithm Engineer at ByteDance.

\bigskip

\noindent\textbf{Zhi Ding (S'88-M'90-SM'95-F'03)} is with the Department of Electrical and Computer Engineering at the University of California, Davis, where he holds the position of distinguished professor. He received his Ph.D. degree in Electrical Engineering from Cornell University in 1990. From 1990 to 2000, he was a faculty member of Auburn University and later, University of Iowa. Prof. Ding joined the College of Engineering at UC Davis in 2000. His major research interests and expertise cover the areas of wireless networking, communications, signal processing, multimedia, and learning. Prof. Ding supervised 40 PhD dissertations since joining UC Davis. His research team of enthusiastic researchers works very closely with industry to solve practical problems and contributes to technological advances in signal processing, wireless networking, and learning. 

Prof. Ding is a Fellow of IEEE and has served as the Chief Information Officer, Chief Marketing Officer, and Parliamentarian of the IEEE Communications Society. He was associate editor for IEEE Transactions on Signal Processing from 1994-1997, 2001-2004, and associate editor of IEEE Signal Processing Letters 2002-2005. He was a member of technical committee on Statistical Signal and Array Processing and member of technical committee on Signal Processing for Communications (1994-2003). Dr. Ding was the General Chair of the 2016 IEEE International Conference on Acoustics, Speech, and Signal Processing and the Technical Program Chair of the 2006 IEEE Globecom. He was also an IEEE {\em Distinguished Lecturer} (Circuits and Systems Society, 2004-06, Communications Society, 2008-09). He served on as IEEE Transactions on Wireless Communications Steering Committee Member (2007-2009) and its Chair (2009-2010). Dr. Ding is a coauthor of the textbook: {\em Modern Digital and Analog Communication Systems}, 5th edition, Oxford University Press, 2019. Prof. Ding received the IEEE Communication Society’s WTC Award in 2012 and the IEEE Communication Society’s Education Award in 2020.

\begin{table}[H]
\centering
\caption{Notations}
\label{tab:notations}
\renewcommand{\arraystretch}{1.02}
\setlength{\tabcolsep}{4pt}
\begin{tabular}{p{0.18\columnwidth}p{0.72\columnwidth}}
\toprule
\textbf{Notation} & \textbf{Description} \\
\midrule
$n$ & Index of a specific channel \\
$N$ & Number of shared channels \\
$M$ & Number of channels chosen at each time \\
$S$ & State of a channel \\
$B$ & Throughput of a channel \\
$S_0,S_1$ & Poor or good channel state \\
$p_{ij}$ & Transition probability between channel states \\
$q$ & Observed CQI level of a channel \\
$K$ & Number of CQI levels \\
$N_{\mathrm{MC}}$ & Number of Monte Carlo simulation runs \\
$p_{i,s}$ & Probability of observing $q=i$ given $S=s$ \\
$a$ & Action taken for a channel \\
$\omega$ & Belief state of a channel \\
\bottomrule
\end{tabular}
\vspace{-0.8em}
\end{table}

\begin{table*}[!p]
\centering
\caption{Terminal average discounted return of every policy in the
Monte Carlo experiments. Variance and standard deviation are calculated across runs, and
the 95\% confidence interval is
$\bar{G}\pm1.96s/\sqrt{N_{\mathrm{MC}}}$. The best performing policy in each
experimental setting, measured by the highest mean return, is
highlighted in bold.}
\label{tab:experiment_statistics}
\scriptsize
\renewcommand{\arraystretch}{0.82}
\setlength{\tabcolsep}{4pt}

\begin{tabular*}{\textwidth}
{@{\extracolsep{\fill}}lrrrr@{}}
\toprule
Policy & Mean & Variance & Std. & 95\% CI \\
\midrule

\multicolumn{5}{@{}l}{
\textbf{Baseline comparison, $\beta=0.5$:}
$N=8$, $M=1$, $K=2$, $T=10$, $N_{\mathrm{MC}}=2000$}\\
VFA             & 0.5366 & 0.0995 & 0.3155 & $[0.5228,0.5504]$ \\
Myopic          & 0.2924 & 0.0776 & 0.2786 & $[0.2802,0.3047]$ \\
Myopic--Rollout & 0.3551 & 0.0692 & 0.2630 & $[0.3435,0.3666]$ \\
AWI0            & 0.2378 & 0.0631 & 0.2512 & $[0.2268,0.2488]$ \\
AWI-0           & 0.3158 & 0.0625 & 0.2499 & $[0.3048,0.3268]$ \\
AWI-1           & 1.2666 & 0.0543 & 0.2329 & $[1.2564,1.2768]$ \\
AWI-2           & 1.3796 & 0.0553 & 0.2352 & $[1.3693,1.3899]$ \\
AWI-3           & \textbf{1.5398} & 0.0608 & 0.2466 & $[1.5290,1.5506]$ \\

\addlinespace[0.5ex]
\multicolumn{5}{@{}l}{
\textbf{Baseline comparison, $\beta=0.9$:}
$N=8$, $M=1$, $K=2$, $T=10$, $N_{\mathrm{MC}}=2000$}\\
VFA             & 3.9476 & 1.8305 & 1.3530 & $[3.8883,4.0069]$ \\
Myopic          & 3.7492 & 1.0499 & 1.0246 & $[3.7043,3.7941]$ \\
Myopic--Rollout & 3.8739 & 1.7837 & 1.3355 & $[3.8154,3.9325]$ \\
AWI0            & 2.5507 & 0.8979 & 0.9476 & $[2.5091,2.5922]$ \\
AWI-0           & 3.1021 & 1.3262 & 1.1516 & $[3.0517,3.1526]$ \\
AWI-1           & 3.8377 & 1.6998 & 1.3038 & $[3.7806,3.8949]$ \\
AWI-2           & 5.0987 & 1.7877 & 1.3370 & $[5.0401,5.1573]$ \\
AWI-3           & \textbf{5.3684} & 1.8659 & 1.3660 & $[5.3086,5.4283]$ \\

\addlinespace[0.5ex]
\multicolumn{5}{@{}l}{
\textbf{Large-system performance, $\beta=0.5$:}
$N=100$, $M=1$, $K=2$, $T=10$, $N_{\mathrm{MC}}=200$}\\
Myopic & 0.3188 & 0.0924 & 0.3039 & $[0.2767,0.3609]$ \\
AWI0   & 0.1971 & 0.0465 & 0.2156 & $[0.1672,0.2270]$ \\
AWI-0  & 0.1732 & 0.0244 & 0.1562 & $[0.1516,0.1949]$ \\
AWI-1  & 0.2686 & 0.0506 & 0.2249 & $[0.2374,0.2997]$ \\
AWI-2  & 0.7476 & 0.0416 & 0.2039 & $[0.7193,0.7758]$ \\
AWI-3  & \textbf{1.5574} & 0.0850 & 0.2915 & $[1.5170,1.5978]$ \\

\addlinespace[0.5ex]
\multicolumn{5}{@{}l}{
\textbf{Large-system performance, $\beta=0.9$:}
$N=100$, $M=1$, $K=2$, $T=10$, $N_{\mathrm{MC}}=200$}\\
Myopic & 3.4302 & 2.1651 & 1.4714 & $[3.2262,3.6341]$ \\
AWI0   & 2.2713 & 1.0263 & 1.0131 & $[2.1309,2.4117]$ \\
AWI-0  & 2.5984 & 1.7427 & 1.3201 & $[2.4154,2.7813]$ \\
AWI-1  & 3.8100 & 1.4334 & 1.1973 & $[3.6441,3.9759]$ \\
AWI-2  & 5.0761 & 1.7275 & 1.3144 & $[4.8939,5.2582]$ \\
AWI-3  & \textbf{5.3530} & 1.7858 & 1.3363 & $[5.1678,5.5382]$ \\

\addlinespace[0.5ex]
\multicolumn{5}{@{}l}{
\textbf{Non-stationary dynamics, $\beta=0.5$:}
$N=8$, $M=1$, $K=2$, $T=10$, $N_{\mathrm{MC}}=500$}\\
Myopic & 0.4158 & 0.0262 & 0.1618 & $[0.4016,0.4300]$ \\
AWI0   & 0.4297 & 0.0069 & 0.0833 & $[0.4224,0.4370]$ \\
AWI-0  & 0.4447 & 0.0066 & 0.0813 & $[0.4376,0.4518]$ \\
AWI-1  & 0.5161 & 0.0044 & 0.0666 & $[0.5103,0.5220]$ \\
AWI-2  & \textbf{0.5575} & 0.0045 & 0.0674 & $[0.5516,0.5634]$ \\
AWI-3  & \textbf{0.5575} & 0.0045 & 0.0674 & $[0.5516,0.5634]$ \\

\addlinespace[0.5ex]
\multicolumn{5}{@{}l}{
\textbf{Non-stationary dynamics, $\beta=0.9$:}
$N=8$, $M=1$, $K=2$, $T=10$, $N_{\mathrm{MC}}=500$}\\
Myopic & 3.2540 & 0.5274 & 0.7262 & $[3.1904,3.3177]$ \\
AWI0   & 3.1249 & 0.7941 & 0.8911 & $[3.0468,3.2030]$ \\
AWI-0  & 3.8700 & 0.6602 & 0.8125 & $[3.7988,3.9412]$ \\
AWI-1  & 4.0234 & 0.6931 & 0.8325 & $[3.9505,4.0964]$ \\
AWI-2  & 4.1303 & 0.7192 & 0.8480 & $[4.0559,4.2046]$ \\
AWI-3  & \textbf{4.1325} & 0.7180 & 0.8473 & $[4.0583,4.2068]$ \\

\addlinespace[0.5ex]
\multicolumn{5}{@{}l}{
\textbf{CQI sensitivity, $K=2$:}
$N=7$, $M=1$, $T=10$, $\beta=0.9$, $N_{\mathrm{MC}}=200$}\\
Myopic & 4.4110 & 0.4342 & 0.6589 & $[4.3197,4.5023]$ \\
AWI0   & 2.9199 & 0.4171 & 0.6458 & $[2.8304,3.0094]$ \\
AWI-0  & 3.4031 & 0.4765 & 0.6903 & $[3.3074,3.4988]$ \\
AWI-1  & 3.6543 & 0.7273 & 0.8528 & $[3.5361,3.7725]$ \\
AWI-2  & 4.1477 & 0.4901 & 0.7001 & $[4.0507,4.2448]$ \\
AWI-3  & \textbf{5.1607} & 0.7168 & 0.8467 & $[5.0433,5.2780]$ \\

\addlinespace[0.5ex]
\multicolumn{5}{@{}l}{
\textbf{CQI sensitivity, $K=3$:}
$N=7$, $M=1$, $T=10$, $\beta=0.9$, $N_{\mathrm{MC}}=200$}\\
Myopic & 4.4098 & 0.4793 & 0.6923 & $[4.3138,4.5057]$ \\
AWI0   & 3.1431 & 0.4925 & 0.7018 & $[3.0458,3.2403]$ \\
AWI-0  & 3.6195 & 0.5143 & 0.7171 & $[3.5201,3.7189]$ \\
AWI-1  & 3.8633 & 0.6296 & 0.7935 & $[3.7533,3.9732]$ \\
AWI-2  & 5.0450 & 0.6110 & 0.7816 & $[4.9367,5.1534]$ \\
AWI-3  & \textbf{5.2136} & 0.6514 & 0.8071 & $[5.1017,5.3254]$ \\

\addlinespace[0.5ex]
\multicolumn{5}{@{}l}{
\textbf{CQI sensitivity, $K=4$:}
$N=7$, $M=1$, $T=10$, $\beta=0.9$, $N_{\mathrm{MC}}=200$}\\
Myopic & 4.4162 & 0.5018 & 0.7084 & $[4.3180,4.5144]$ \\
AWI0   & 2.8961 & 0.4387 & 0.6624 & $[2.8043,2.9879]$ \\
AWI-0  & 3.6211 & 0.4605 & 0.6786 & $[3.5270,3.7151]$ \\
AWI-1  & 3.8232 & 0.5911 & 0.7688 & $[3.7166,3.9297]$ \\
AWI-2  & 5.1443 & 0.6131 & 0.7830 & $[5.0358,5.2528]$ \\
AWI-3  & \textbf{5.2361} & 0.5988 & 0.7738 & $[5.1289,5.3434]$ \\

\addlinespace[0.5ex]
\multicolumn{5}{@{}l}{
\textbf{CQI sensitivity, $K=5$:}
$N=7$, $M=1$, $T=10$, $\beta=0.9$, $N_{\mathrm{MC}}=200$}\\
Myopic & 4.3912 & 0.5181 & 0.7198 & $[4.2914,4.4909]$ \\
AWI0   & 2.8576 & 0.5651 & 0.7517 & $[2.7534,2.9617]$ \\
AWI-0  & 3.6667 & 0.4599 & 0.6782 & $[3.5728,3.7607]$ \\
AWI-1  & 3.9382 & 0.6181 & 0.7862 & $[3.8293,4.0472]$ \\
AWI-2  & 4.9175 & 0.6392 & 0.7995 & $[4.8067,5.0283]$ \\
AWI-3  & \textbf{5.2039} & 0.6600 & 0.8124 & $[5.0913,5.3165]$ \\

\addlinespace[0.5ex]
\multicolumn{5}{@{}l}{
\textbf{Extreme discount factor:}
$N=8$, $M=1$, $K=2$, $T=10$, $\beta=1$, $N_{\mathrm{MC}}=2000$}\\
Myopic & 6.1310 & 2.2764 & 1.5088 & $[6.0648,6.1971]$ \\
AWI0   & 3.9563 & 2.4173 & 1.5548 & $[3.8882,4.0245]$ \\
AWI-0  & 2.8960 & 5.4950 & 2.3441 & $[2.7932,2.9987]$ \\
AWI-1  & 3.8021 & 4.9042 & 2.2145 & $[3.7051,3.8992]$ \\
AWI-2  & 4.1842 & 4.6171 & 2.1487 & $[4.0900,4.2784]$ \\
AWI-3  & \textbf{8.0705} & 5.4734 & 2.3395 & $[7.9680,8.1731]$ \\

\bottomrule
\end{tabular*}
\end{table*}

\begin{table*}[!t]
\centering
\caption{Threshold optimality and indexability classification of the
experimental instances. TO means threshold optimality and which is marked ``Yes'' only when every arm
satisfies the condition.}
\label{tab:threshold_indexability_cases}

\resizebox{\textwidth}{!}{
\begin{tabular}{lcccl}
\hline
Experiment
& $\beta$
& TO
& Indexability
& Classification \\
\hline

Baseline comparison
& 0.5 & No & No
& Neither condition \\

Baseline comparison
& 0.9 & No & No
& Neither condition \\

Exact DP benchmark
& 0.9 & No & No
& Neither condition \\

Large system
& 0.5 & No & No
& Neither condition \\

Large system
& 0.9 & No & No
& Neither condition \\

Non-stationary
& 0.5 & Yes & Yes
& Both conditions \\

Non-stationary
& 0.9 & No & No
& Neither condition \\

CQI sensitivity, $K=2$
& 0.9 & Yes & No
& Threshold optimality only \\

CQI sensitivity, $K=3$
& 0.9 & Yes & No
& Threshold optimality only \\

CQI sensitivity, $K=4$
& 0.9 & Yes & No
& Threshold optimality only \\

CQI sensitivity, $K=5$
& 0.9 & Yes & No
& Threshold optimality only \\

Extreme discount factor
& 1 & No & No
& Neither condition \\

\hline
\end{tabular}
}
\end{table*}

\begin{table*}[!t]
\centering
\caption{Comparison with the binary-feedback RMAB in~\cite{liu_low-complexity_2022} and preliminary conference paper~\cite{elmaghraby_femtocell_2018}.}
\label{tab:comparison}
\small
\renewcommand{\arraystretch}{1.15}
\begin{tabularx}{\textwidth}{@{}p{0.21\textwidth}XXX@{}}
\toprule
Feature
& Binary-feedback RMAB~\cite{liu_low-complexity_2022}
& Preliminary conference~\cite{elmaghraby_femtocell_2018}
& This work \\
\midrule

CQI Observation Model
& Binary feedback
& Multi-level CQI
& Multi-level CQI \\
\midrule

Belief Transition Structure
& Bayesian update with binary feedback
& Bayesian update with multi-level CQI
& Bayesian update with multi-level CQI \\
\midrule

Derivative bounds
& No convergence analysis of approximation errors
& \textup{N/A}
& Approximation errors converge at a geometric rate \\
\midrule

Threshold/Indexability Conditions
& Yes
& \textup{N/A}
& Yes \\
\midrule

Iterative AWI computation
& No complexity analysis
& Partly, $n_{\mathrm{iter}}=0$ and $1$; no complexity analysis
& Complete analysis \\
\midrule

Simulation extensions
& Only Myopic baseline
& No baseline comparison
& Stronger and multiple baselines \\

\bottomrule
\end{tabularx}
\end{table*}

\begin{table*}[t]
    \centering
    \caption{Parameters of the two baseline systems. System-1 corresponds to $\beta=0.5$, while System-2 corresponds to $\beta=0.9$.}
    \label{tab:baseline_parameters}
    
    \small
    \renewcommand{\arraystretch}{1.3} 
    
    \begin{tabularx}{\textwidth}{l c X}
        \toprule
        \textbf{System Model} & \textbf{Parameter} & \textbf{Value Vector} \\
        \midrule
        
        \textbf{System-1} 
        & $\left\{p_{11}^{(n)}\right\}_{n=1}^{8}$ 
        & $\{0.91,\, 0.20,\, 0.47,\, 0.48,\, 0.60,\, 0.34,\, 0.69,\, 0.40\}$ \\
        \midrule
        $(\beta=0.5)$ 
        & $\left\{p_{01}^{(n)}\right\}_{n=1}^{8}$ 
        & $\{0.09,\, 0.87,\, 0.41,\, 0.11,\, 0.26,\, 0.28,\, 0.33,\, 0.47\}$ \\
        \midrule
        & $\left\{B_n\right\}_{n=1}^{8}$ 
        & $\{0.96,\, 0.78,\, 1.00,\, 0.83,\, 0.79,\, 0.81,\, 0.89,\, 0.88\}$ \\
        \midrule
        & $\left\{p_{i,1}^n\right\}_{n=1}^{8}$ 
        & $\{0.77,\, 0.69,\, 0.66,\, 0.76,\, 0.82,\, 0.74,\, 0.71,\, 0.76\}$ \\
        \midrule
        
        \textbf{System-2} 
        & $\left\{p_{11}^{(n)}\right\}_{n=1}^{8}$ 
        & $\{0.95,\, 0.72,\, 0.43,\, 0.79,\, 0.71,\, 0.49,\, 0.22,\, 0.97\}$ \\
        \midrule
        $(\beta=0.9)$ 
        & $\left\{p_{01}^{(n)}\right\}_{n=1}^{8}$ 
        & $\{0.11,\, 0.89,\, 0.55,\, 0.08,\, 0.12,\, 0.62,\, 0.52,\, 0.03\}$ \\
        \midrule
        & $\left\{B_n\right\}_{n=1}^{8}$ 
        & $\{0.98,\, 0.83,\, 0.88,\, 1.00,\, 0.87,\, 0.82,\, 0.84,\, 0.92\}$ \\
        \midrule
        & $\left\{p_{i,1}^n\right\}_{n=1}^{8}$ 
        & $\{0.76,\, 0.76,\, 0.71,\, 0.73,\, 0.58,\, 0.78,\, 0.75,\, 0.60\}$ \\
        
        \bottomrule
    \end{tabularx}

\end{table*}


 





\end{document}